\documentclass[sort&compress,sn-mathphys,Numbered]{sn-jnl}

\usepackage{multirow}%
\usepackage{amsthm}%
\usepackage{amsmath}%
\usepackage{mathrsfs}%
\usepackage[title]{appendix}%
\usepackage{xcolor}%
\usepackage{textcomp}%
\usepackage{manyfoot}%
\usepackage{booktabs}%
\usepackage{algorithm}%
\usepackage{algorithmicx}%
\usepackage{algpseudocode}%
\usepackage{listings}%
\usepackage{amsmath,amssymb,amsfonts,amscd,amsbsy,amsxtra}
\usepackage{booktabs} 
\usepackage{caption} 
\usepackage{subcaption} 
\usepackage{array}
\usepackage{pgfplots}
\usepackage[all]{nowidow}
\usepackage[utf8]{inputenc}
\usepackage{multicol}
\usepackage{algpseudocode,algorithm,algorithmicx}
\usepackage{hyperref}
\usepackage{cleveref}
\usepackage[inline]{enumitem} 
\usepackage{csquotes}
\usepackage{mathtools}
\usepackage{bm}
\usepackage{bbm}
\usepackage[mathscr]{euscript} 
\usepackage{doc}
\usepackage{url}
\usepackage{soul}

\usepackage{paralist}

\usepackage{cancel}

\usepackage[textsize=tiny]{todonotes}

\definecolor{mediumred}{RGB}{204,0,0}
\newcommand\rev[1]{{\textcolor{black}{#1}}}

\usepackage{tikz}
\usetikzlibrary{er,positioning,bayesnet}
\pgfplotsset{compat=1.16}
\usepackage{graphicx}
\definecolor{blue}{HTML}{1F77B4}
\definecolor{orange}{HTML}{FF7F0E}
\definecolor{green}{HTML}{2CA02C}

\allowdisplaybreaks[4]

\numberwithin{equation}{section}
\newtheorem{thm}{theorem}[section]
\numberwithin{thm}{section}
\newtheorem{lemma}[thm]{Lemma}
\newtheorem{theorem}[thm]{Theorem}
\newtheorem{corollary}[thm]{Corollary}
\newtheorem{proposition}[thm]{Proposition}
\newtheorem{remark}[thm]{Remark}
\newtheorem{example}[thm]{Example}
\newtheorem{definition}[thm]{Definition}

\DeclareMathOperator{\dist}{dist}
\DeclareMathOperator{\Span}{span}

\newcommand\xqed[1]{%
  \leavevmode\unskip\penalty9999 \hbox{}\nobreak\hfill
  \quad\hbox{#1}}
\newcommand\demodrei{\xqed{\small $\triangle$}}

\newcommand\HWS{\mathrm{hws}}
\newcommand\even{\mathrm{evn}}
\newcommand\odd{\mathrm{odd}}
\newcommand\eo{\mathrm{eo}}

\newcommand\Wke{W_k^\even}
\newcommand\Wko{W_k^\odd}

\newcommand{\cB}{\mathcal{B}}
\newcommand{\cK}{\mathcal{K}}

\newcommand{\cO}{\mathcal{O}}
\newcommand{\cP}{\mathcal{P}}
\newcommand{\cQ}{\mathcal{Q}}
\newcommand{\cU}{\mathcal{U}}

\newcommand{\cW}{\mathcal{W}}

\newcommand{\OP}{\Omega_{\cP}}

\newcommand{\C}{\mathbb{C}}
\newcommand{\N}{\mathbb{N}}
\newcommand{\R}{\mathbb{R}}
\newcommand{\Z}{\mathbb{Z}}

\newcommand{\bc}{\bm{c}}

\newcommand{\bx}{\bm{x}}

\newcommand{\bU}{\bm{U}}
\newcommand{\bV}{\bm{V}}
\newcommand{\bX}{\bm{X}}

\newcommand{\bmu}{\bm{\mu}}

\newcommand{\bPsi}{\bm{\Psi}}
\newcommand{\bSigma}{\bm{\Sigma}}

\usepackage{MnSymbol}

\begin{document}
\title[Kolmogorov $N$-with for the transport problem]{The Kolmogorov $N$-width for linear transport: Exact representation and the influence of the data}
%
\author[1]{\fnm{Florian} \sur{Arbes}}
 \email{florian.arbes@ife.no}
\author[2]{\fnm{Constantin} \sur{Greif}}\email{constantin.greif@uni-ulm.de}
\author[2]{\fnm{Karsten} \sur{Urban}}\email{karsten.urban@uni-ulm.de}
\affil[1]{\orgdiv{Computational Materials Processing}, \orgname{IFE, Institute for Energy Technology}, \orgaddress{\street{Instituttveien 18}, \city{Kjeller}, \postcode{2007}, \country{Norway}}}
\affil[2]{\orgdiv{Institute of Numerical Mathematics}, \orgname{Ulm University}, \orgaddress{\street{Helmholtzstr. 20}, \city{Ulm}, \postcode{89081}, \country{Germany}}}

\abstract{
The Kolmogorov $N$-width describes the best possible error one can achieve by elements of an $N$-dimensional linear space.
Its decay has extensively been studied in Approximation Theory and for the solution of Partial Differential Equations (PDEs).
Particular interest has occurred within Model Order Reduction (MOR) of parameterized PDEs e.g.\ by the Reduced Basis Method (RBM).

While it is known that the $N$-width decays exponentially fast (and thus admits efficient MOR) for certain problems, there are examples of the linear transport and the wave equation, where the decay rate deteriorates to $N^{-1/2}$.
On the other hand, it is widely accepted that a smooth parameter dependence admits a fast decay of the $N$-width.
However, a detailed analysis of the influence of properties of the data (such as regularity or slope) on the rate of the $N$-width seems to lack.

In this paper, \rev{we state that the optimal linear space is a direct sum of shift-isometric eigenspaces corresponding to the largest eigenvalues, yielding an exact representation of the $N$-width as their sum. For the linear transport problem, which is modeled by half-wave symmetric initial and boundary conditions $g$, we obtain such an optimal decomposition by sorted trigonometric functions with eigenvalues that  match the Fourier coefficients of $g$.
Further the sorted eigenfunctions give for normalized $g\in H^{r}$ of broken order $r>0$ the sharp upper bound of the $N$-width, which is a reciprocal of a certain power sum. Yet for ease, we also provide the decay $(\pi N)^{-r}$, obtained by the non-optimal space of ordering the trigonometric functions by frequency rather than by eigenvalue.}

Our theoretical investigations are complemented by numerical experiments which confirm the sharpness of our bounds and give additional quantitative insight.}
\maketitle              
\section{Introduction} \label{introduction}
The Kolmogorov $N$-width describes the best possible error one can achieve by a linear approximation with $N\in\N$ degrees of freedom, i.e.\ by elements of the best possible $N$-dimensional linear space \cite{kolmogoroffUberBesteAnnaherung1936}. The arising optimal space in the sense of Kolmogorov can often not explicitly be constructed, at least not in a reasonable (computing) time. On the other hand, however, the decay rate of the $N$-width tells us if a given set can be well-approximated by a linear method, or not. This is a classical question in Approximation Theory and has been widely studied in the literature, see e.g.\ \cite{pinkusNWidthsApproximationTheory1985,devore1998nonlinear,floaterBestLowrankApproximations2021a,bressanBestConstantsL22021}, this list being far from complete.

Particular interest has been devoted to the case when the set to be approximated is given by solutions of certain equations, e.g.\ Partial Differential Equations (PDEs) with different data (one might think of the domain, coefficients, right-hand side loadings, initial- and/or boundary conditions), which might be considered as parameters \cite{cohen2010convergence,cohenKolmogorovWidthsHolomorphic2015,MELENK2000272}. 
In that direction, model order reduction of parametric PDEs (PPDE) has become a field of very intensive research, also with many very relevant real-life applications \cite{bennerModelReductionApproximation2017,hesthavenCertifiedReducedBasis2015,quarteroniReducedBasisMethods2015}.
A prominent example is the Reduced Basis Method (RBM), where a PPDE is aimed to be reduced to an $N$-dimensional linear space in order to allow multi-query (w.r.t.\ different parameter values) and/or realtime (embedded systems, cloud computing) applications.
The reduced $N$-dimensional system is determined in an offline training phase using sufficiently accurate detailed numerical solutions by any standard method. 
In this framework, the question arises, if a given PPDE can be well-reduced by means of the RBM or not.
Since it has been proven in \cite{binevConvergenceRatesGreedy2011} that the offline reduced basis generation using a Greedy method realizes the same asymptotic rate of decay as the Kolmogorov $N$-width, one is left with the investigation of the decay for PPDEs to decide whether the RBM is suitable for a given PPDE, or not.
Also in that direction, there is a significant amount of literature, e.g.\ \cite{bachmayrKolmogorovWidthsLowrank2015a,devoreChapterTheoreticalFoundation2017,maday2002priori,maday2006reduced,MR2430350,MR2877366,Lassila2013,ohlbergerReducedBasisMethods2016c,MPT,MADAY2002289}, just to name a few.
Roughly speaking, it was shown there that a PPDE admits a fast decay of the $N$-width if the solution depends smoothly on the parameter, which is, e.g.\ known for elliptic and parabolic problems which allow for a separation of the parameters from the physical variables.
As a rule of thumb: \enquote{holomorphic dependence admits exponential decay}.

However, when leaving the \emph{nice} realm of such PPDEs, the situation becomes dramatically worse.
It has for example been shown that the decay may drop down to $N^{-s}$, $0<s<1$, for the linear transport equation \cite{ohlbergerReducedBasisMethods2016c} and the wave equation \cite{greifDecayKolmogorovNwidth2019}.
However, the problems considered in the latter papers are quite specific examples yielding to a non-smooth dependence of the solution in terms of the parameter (the velocity in \cite{ohlbergerReducedBasisMethods2016c} and the wave speed in \cite{greifDecayKolmogorovNwidth2019}). It was also demonstrated that the decay not only depends on the PDE, but also on the underlying physics, e.g.\ alloy compositions in case of solidification problems \cite{arbesModelOrderReduction2022a}. For problems of such type (transport, transport-dominated, hyperbolic), the above quoted rule of thumb remains true.

This is why we are interested in the \emph{exact} dependence of the decay rate of the $N$-width in terms of the data / parameters of the problem.
To our own surprise we could not find corresponding results in the literature.
In \cite{MR2877366,MPT,MADAY2002289}, the fast decay is shown using techniques from interpolation proving that a Greedy-type selection selects the optimal nodes.
The positive result in \cite[Thm.\ 3.1]{ohlbergerReducedBasisMethods2016c} has been deduced by using the decay of the complex power series.

We consider the linear transport problem whose solution is given by the characteristics in terms of initial and boundary conditions.
Hence, we can reduce ourselves to approximate the mapping $x\mapsto g(x-\mu)$, where $\mu$ is the parameter and $x\in\Omega$, which is the domain on which the PPDE is posed; $g$ is the real-valued univariate function modeling initial and boundary conditions.
To this end, we use the Fourier series approximation, which allows us to incorporate the parametric shift by $\mu$ into the approximation spaces.
We derive exact representations of the $N$-width for certain classes (half-wave symmetric -HWS- functions) and give estimates in terms of the regularity and the slope of the function $g$.

In Section~\ref{Sec:NWidth}, we use Def.~\ref{ortho-subspace-basis} to construct approximation spaces with which we can derive exact representations for the $N$-width for HWS functions and sharp estimates in the general case.

This paper is organized as follows. 
\rev{In Section~\ref{Sec:Preliminaries}, we gather preliminaries on the linear transport equation, the \( N \)-width, and some facts from Fourier Analysis, which we will use in the course of the paper. Section 2 sets the stage for the various sections, which gradually become more specific.
\\
In Section~\ref{section3}, we prove that for parametric problems, the \( L_2 \)-optimal spaces (mean-squared error in the parameter) consist of spectral spaces, which, to achieve \( L_{\infty} \)-optimality (worst-case in the parameter), must also be shift-isometric. The \( N \)-width is then simply a sum of all eigenvalues minus the \( N \) largest ones. Here, the connection to Proper Orthogonal
Decomposition (POD) should be mentioned.
\\
In Section~\ref{Sec:NWidth}, we specifically address the linear transport problem (cf. Preliminaries 2.1) with half-wave symmetric initial and boundary values (cf. Preliminaries 2.3) and prove the aforementioned properties. It turns out that, in this setting, the eigenvalues correspond to the Fourier coefficients of the initial and boundary values.
\\
In Section~\ref{Sec:Smoothness}, we further specifically involve fractional Sobolev regularity (cf. Preliminaries 2.4). This results in a sharp upper bound of the \( N \)-width, which is a reciprocal of a certain power sum scaled by the \( H^r \)-seminorm. Yet for ease, we provide the decay \( (\pi N)^{-r} \), obtained by the non-optimal space of ordering the trigonometric functions by frequency rather than eigenvalue.\\
For specific data, we can further define the exact rate, which we elaborate in Section 5.2. As part of this, we obtain a decay of $c_m N^{-m-1/2}, m \in \N$, for piecewise polynomial data $g \in H^m$,
such that $g^{(m)}$ represents a staircase function with regularity less than $1/2$. We conclude the theoretical section with a discussion of what happens as regularity \( r \to \infty \).\\
Some results of our numerical investigations are presented in Section~\ref{Sec:experiments}.
First we visualize our results for the exact form and the asymptotic rate of the error bound of a jump discontinuity. 
Then, we quantify the decay constants for a set functions with known error bounds that have approximately the same slope but a varying degree of smoothness.
Analogously we fix the regularity and vary the slope in a third experiment.
Lastly, in Section~\ref{sec:beyhond_symm} we consider almost arbitrary functions to investigate how the results extend beyond the symmetry constraints.
}
The paper finishes with some conclusions in Section~\ref{Sec:Outlook}.

\rev{The four main results in the paper are first, \underline{Corollary~\ref{cor-Nwidths-equal}} which gives for a linear $N$-dimensional shift-isometric space an exact representation of both $N$-widths. Second, \underline{Theorem~\ref{Theorem:Optimal-simple}} which states that we get an exact representation of the $N$-width for the linear transport equation with half-wave symmetric data. Third, \underline{Theorem~\ref{sol-manifold-decay}} which gives the sharp upper bound of the $N$-width for the linear transport problem for data $g \in H^r$.
And finally, the \underline{numerical experiments} show that the decay constants heavily depend on the slope of the function and demonstrate how small changes can increase regularity and thus lead to much faster error decays.}

\section{Preliminaries}\label{Sec:Preliminaries}
\subsection{The linear transport equation} \label{section-transport}
We consider the univariate linear transport equation with velocity $\mu$, which is interpreted as a parameter, i.e., we seek a function $\Phi(\cdot,\cdot;\mu): I\times\Omega\to\mathbb{R}$ such that\footnote{We restrict ourselves to the homogeneous case for simplicity.
    One could also consider a right-hand side $f(t,x;\mu)\not\equiv 0$, which would also impact the rate of approximation of the solution.}%
\begin{subequations}
    \label{eq:transport}
    \begin{align}
        \partial_{t} \Phi(t,x;\mu)
        + \mu \, \partial_x \Phi(t,x;\mu) & = 0,
                                       &                 (t,x) & \in I \times \Omega,
        \label{eq:transport:1}                                                       \\
        \Phi(0,x;\mu)                     & = g(x),
                                       &                 x & \in \Omega,
        \label{eq:transport:2}                                                       \\
        \Phi(t,0;\mu)                     & = g(-\mu\, t),
                                       &                 t & \in I,
        \label{eq:transport:3}
    \end{align}
\end{subequations}
where $I:=(0,1)$ is the time interval and $\Omega:=(0,1)$ the spatial domain.\footnote{Our analysis is restricted to the 1D-case, but some results can be extended to higher dimensions.} The velocity can be chosen in a compact interval $\mu\in\cP:=[0,1]$.
The initial and boundary conditions, \eqref{eq:transport:2} and \eqref{eq:transport:3}, respectively, are given in terms of a function $g:\OP\to\R$, whose properties will be relevant in the sequel.
Here, 
$$\OP:=\{x-\mu t: x \in \Omega, \mu \in \cP, t \in I \}=(-1,1)$$ 
is the domain on which $g$ needs to be defined in order to obtain a well-posed problem \eqref{eq:transport} for every parameter.
The solution of \eqref{eq:transport} is well-known to read $\Phi(t,x;\mu)=g(x-\mu\, t),$ $(t,x) \in I \times \Omega$.
We are particularly interested in the solution at the final time $t=1$, i.e., 
\vspace*{-0.3cm}
\begin{align}\label{eq:functiontobeapproximated}
    u_\mu(x):= \Phi(1,x;\mu) = g(x-\mu),
    \qquad x\in\Omega,
\end{align} 
and consider the low regularity case, i.e., we only assume that $g \in L_2(\OP)=\rev{L_2(\OP;\R)}$, and therefore $u_\mu\in L_2(\Omega)=\rev{L_2(\Omega;\R)}$, without additional smoothness.

\begin{remark}
    Often, the time $t \in I$ is also seen as a parameter.
    But our considerations are not restricted to the final time, since then for a given $t \in I$ and $\mu \in \cP$, we can define the new parameter $\tilde{\mu} := t\,  \mu \in [0,1)$ and get $u_{\tilde{\mu}}= g(\cdot- \tilde{\mu})=\Phi(t,\cdot;\mu)$.\hfill$\diamond$
\end{remark}

\subsection{Linear approximation: The \textit{N}-width}
The specific focus of this paper is the approximation rate provided by linear subspaces.
In particular, we are considering $N$-dimensional subspaces which are \enquote{optimal} to approximate $u_\mu$ for all parameters $\mu\in\cP$ in an appropriate manner.  Our aim is to study the dependency of the rate of approximation w.r.t.\ the data of the problem, namely initial and boundary conditions modeled by the function $g$. The parameter set $\cP$ is fixed. Hence, we aim at approximating the \enquote{solution manifold}
\begin{align}\label{eq:SolManifold}
	\cU^g := \{ u_\mu=g(\cdot-\mu):\, \mu\in\cP\}
	\subset L_2(\Omega).
\end{align}
The maybe most classical setting is the worst case scenario w.r.t.\ the parameter yielding the classical \emph{Kolmogorov $N$-width} \cite{kolmogoroffUberBesteAnnaherung1936} defined as 
\begin{align}\label{eq:KolNWdth}
    d_N(\cU^g)
      :=& \kern-4pt\inf_{\substack{V_N\subset L_2(\Omega) \\ \dim(V_N)=N}}
    \sup_{\mu\in\cP}
    \inf_{\tilde v_N\in V_N}
    \| u_\mu - \tilde v_N\|_{L_2(\Omega)}
    \\
    =& \kern-4pt\inf_{\substack{V_N\subset L_2(\Omega)     \\ \dim(V_N)=N}}
    \dist(V_N, \cU^g)_{L_\infty(\cP;L_2(\Omega))}.
    \nonumber
\end{align}
The dependence on $g$ will be crucial below.

\begin{remark} There are several results concerning the decay of $d_N( \cU^g )$ for the linear transport problem \eqref{eq:transport}.\vspace*{-0.05\baselineskip}
    \begin{compactenum}[(i)]
        \item In \cite{ohlbergerReducedBasisMethods2016c} it was shown that $d_N( \cU^g )$ decays as $N^{-1/2}$, i.e., very slowly for the specific choice $g=\chi_{[0,1]}$, namely for initial and boundary conditions involving a jump.
        \item On the other hand, one can achieve exponential decay, i.e., $d_N( \cU^g ) \lesssim e^{-\alpha N}$ for some $\alpha>0$ if the function $g$ is analytic.
              This can be seen by considering a truncated power series in the complex plane \cite[Thm.\ 3.1]{ohlbergerReducedBasisMethods2016c}.\hfill$\diamond$
    \end{compactenum}
\end{remark}
$\text{ }$\\
Our main focus in this paper is to study the decay of the $N$-width w.r.t.\ properties of the function $g$, in particular we want to detail the influence of the regularity of $g$ on the decay of the $N$-width.
In addition to the \enquote{worst-case in the parameter} Kolmogorov $N$-width $d_N( \cU^g )$, which measures the error in $L_\infty(\cP;L_2(\Omega))$, we will also consider \enquote{mean-squared error in the parameter}, i.e., $L_2(\cP;L_2(\Omega))$ w.r.t.\ a probability measure, i.e., $\int_{\cP} \rho(\mu) d \mu=1$, \rev{with density function $\rho(\mu) \equiv \tfrac{1}{\lvert \cP \rvert}$}, which we call \emph{$L_2$-average $N$-width} defined as
\begin{align}\label{eq:L2NWdth}
    \delta_N( \cU^g )
      :=& \inf_{\substack{V_N\subset L_2(\Omega) \\ \dim(V_N)=N}} \!
    \Big\{ \rev{\tfrac{1}{\lvert \cP \rvert}}
    \int_{\cP}
    \inf_{\tilde v_N\in V_N}
    \| u_\mu - \tilde v_N\|_{L_2(\Omega)}^2
    \, d\mu
    \Big\}^{1/2}
    \\
     =& \inf_{\substack{V_N\subset L_2(\Omega)  \\ \dim(V_N)=N}}
    \dist(V_N, \cU^g )_{L_2(\cP;L_2(\Omega))}.
    \nonumber
\end{align}

\begin{remark}\label{Remark:Widths}
    For later reference, we collect some properties of the $N$-widths.
    \begin{compactenum}[(i)]
        \item Let $P_{N}: L_2(\Omega)\to V_N$ denote the orthogonal projection onto $V_N$. Then,
              \begin{subequations}\label{dN_2}
                  \begin{align}
                      d_N( \cU^g ) 
                      &= \inf_{{\substack{V_N\subset L_2(\Omega)             \\ \operatorname{dim} (V_N)=N}}}
                      \| u_\mu - P_{N} u_\mu \|_{L_{\infty}(\cP;L_2(\Omega))},\\
                      \delta_N( \cU^g ) 
                      &= \inf_{{\substack{V_N\subset L_2(\Omega)    \\ \dim(V_N)=N}}}
                      \| u_\mu - P_{N}u_\mu \|_{L_{2}(\cP;L_2(\Omega))}.\label{deltaN_2}
                  \end{align}
              \end{subequations}
        \item As $\| w\|_{L_2(\cP)} \le  \| w\|_{L_\infty(\cP)}$ for $w\in L_\infty(\cP)$, we get $\delta_N(\cU^g)\le d_N(\cU^g)$.
        \rev{
        \item For $\alpha\in\R$, it holds that $d_N(\cU^{\alpha g})=|\alpha|\, d_N(\cU^{g})$ and $\delta_N(\cU^{\alpha g})=|\alpha|\, \delta_N(\cU^{g})$.
        \item Let $X(\OP)\subset L_2(\OP)$ be some subspace with squared norm $\|\cdot\|_{X(\OP)}^2 = \lvert \, \cdot \, \rvert_{X(\OP)}^2 + \|\cdot\|_{L_2(\OP)}^2 $ in the sense that $X(\OP)=\{ v\in L_2(\OP):\, \lvert \, v \, \rvert_{X(\OP)}<\infty\}$, and define 
        \begin{align}\label{eq:def:ball}
        		B_{\varrho} = B_{\varrho}[X(\OP)] := \{ v\in X(\OP):\, \lvert \, v \, \rvert_{X(\OP)} \le\varrho\},
        \end{align}
        then $\qquad \sup_{g\in B_{\varrho}} \frac{d_N(\cU^g)}{\varrho} =  \sup_{g\in X(\OP)} \frac{d_N(\cU^g)}{\lvert g \rvert_{X(\OP)}}$.
        }
    \end{compactenum}
\rev{
\begin{proof}
	We shall prove the last item. First, it is obvious that $\sup_{g\in B_{\varrho}} d_N(\cU^g) \ge \sup_{g\in \partial B_{\varrho}} d_N(\cU^g)$. Next, by (iii)
	\begin{align*}
		\sup_{g\in \partial B_{\varrho}} d_N(\cU^g)
		&= \sup_{g\in X} \, d_N \big( \cU^{\varrho \tfrac{g}{\lvert g\rvert_{X(\OP)}}} \big)
		= \sup_{g\in X} \tfrac{\varrho}{\lvert g\rvert_{X(\OP)}} d_N(\cU^g)
		\ge \sup_{g\in B_{\varrho}} d_N(\cU^g),
	\end{align*}
	i.e., $\sup_{g\in B_{\varrho}} d_N(\cU^g) = \sup_{g\in \partial B_{\varrho}} d_N(\cU^g)$, which proves the claim.
\end{proof}}
\end{remark}

\subsection{Fourier Analysis and \rev{half-wave symmetry}}\label{Subsec:Fourier}
Our major tool for determining the decay of the $N$-widths is Fourier Analysis. Hence, we shall always assume that $g$ is periodic on the larger domain $\OP$, which is no restriction for the transport problem under consideration. We collect the main ingredients needed for the sequel of this paper.  Recall that for the above model problem, we have $I=(0,1)$, $\Omega=(0,1)$, $\cP=[0,1]$ and $\OP=(-1,1)$, but the analysis is not restricted to that case. We shall use the Fourier series of $L_2(\OP)$-functions, namely \rev{
\begin{align}\label{eq:Fourier}
	g =  \sum_{k \in \Z} \hat{c}_{k}(g) \, \psi_k ,
	\qquad 
	\psi_k:=\tfrac{1}{\sqrt{2}} {\textstyle e}^{i \pi k \cdot},
\end{align}
where the Fourier coefficients are known as $\hat{c}_{0}(g)=\tfrac{1}{\sqrt{2}} \int_{\OP} g(x)\, dx$ as well as for $k \in \N$ 
\begin{align}\label{eq:ckg}
	\hat{c}_{k}(g)  
	:=  \langle g, \psi_k \rangle_{L_2(\OP;\C)}  
	\qquad 
	\text{and}\quad 
	\hat{c}_{-k}(g) = \overline{\hat{c}_{k}(g)}, 
\end{align}
since $g$ is real-valued. Consequently, the set $\big\{ \frac{\psi_k - \psi_{-k}}{\footnotesize{\sqrt{2}} i}, \frac{\psi_k + \psi_{-k}}{\footnotesize{\sqrt{2}}} \big\}_{k\in \N_0} = \big\{ \sin(k\pi \cdot), \cos(k\pi \cdot) \big\}_{k\in \N_0}$ is an ONB for $L_2(\OP)$.}%
We consider the space $L_2(\OP)$ corresponding to signals of wave-length $2$. Thus, the half-wave length is $1$, which is used in the following definition, whose notion is well-known in electrical engineering (see e.g.\ \cite{attenboroughMathematicsElectricalEngineering2003}) and turns out to be crucial for the subsequent analysis.
\smallskip

\begin{definition}\label{Def:HWS}
     We call $g\in L_2(\OP)$ \emph{even half-wave symmetric (even HWS, $g\in L_2^\even$)}, if $g(x) = g(x + 1)$ for almost all $x \in [-1,0]$, and \emph{odd half-wave symmetric (odd HWS, $g\in L_2^\odd$)}, if $g(x) = -g(x + 1)$ for almost all $x \in [-1,0]$. A function is called \emph{half-wave symmetric (HWS, $g\in L_2^\HWS$)}, if it is either even or odd HWS.\hfill$\diamond$
 \end{definition}
 \smallskip
 
Using the Fourier expansion, it can readily be seen that any $g\in L_2(\OP)$ can be decomposed into an even HWS and and odd HWS part, i.e., $g=g^\even+g^\odd$, where $g^\even\in L_2^\even$ and $g^\odd\in L_2^\odd$, as functions on $\OP$, admit the Fourier expansion \rev{\eqref{eq:Fourier-hws}.} 

\rev{
\begin{remark}\label{RemFourierhws}
	The Fourier expansion simplifies for HWS functions. In fact, 
	\begin{align}
	\label{eq:Fourier-hws}
		g^\even=\sum_{k\in\Z} \hat{c}_{2k}(g^\even) \, \psi_{2k} \in L_2^\even, \,\quad\,
		g^\odd=\sum_{k\in\Z} \hat{c}_{2k-1}(g^\odd) \, \psi_{2k-1}  \in L_2^\odd ,
	\end{align}
	and for $g=g^\even+g^\odd$ we have $\hat{c}_{2k}(g) =\hat{c}_{2k}(g^\even), \, \hat{c}_{2k-1}(g) = \hat{c}_{2k-1}(g^\odd),  \, k\in\Z$.\newline
	This can also be expressed as follows: $g\in L_2^\even$ if and only if $\hat{c}_{2k-1}(g)=0$ and $g\in L_2^\odd$ if and only if $\hat{c}_{2k}(g)=0$ for all $k\in\Z$.
	\hfill$\diamond$
\end{remark}
}

We shall use this decomposition in order to determine the decay of the Kolmogorov $N$-width by splitting $g$ into its even HWS and odd HWS part and then estimating the $N$-width for both of these parts. For later reference, we collect the facts
\begin{align}\label{eq:L2hws}
    L_2(\OP) = L_2^\even \oplus L_2^\odd,
    \qquad
    L_2^\HWS = L_2^\even \, \cupdot \, L_2^\odd,
    \qquad
    	L_2^\HWS \subsetneq L_2(\OP).
\end{align}
\rev{
\begin{remark}\label{Rem:L2Omega}
\begin{compactenum}[(a)]
	\item The sets $\Psi^\even:= \big\{ \sin(2k\pi \cdot), \cos(2k\pi \cdot) \big\}_{k\in \N_0}$ and $\Psi^\odd:= \big\{ \sin((2k-1)\pi \cdot), \cos((2k-1)\pi \cdot) \big\}_{k\in \N}$ are orthonormal bases (ONB) for $L_2^\even$ and $L_2^\odd$ (i.e., on $\OP$), respectively. 
	\item We shall also need orthonormal basis functions in $L_2(\Omega)$ and set $\tilde\psi^\even_0 := 1$ and for $k \in \N$,
	\begin{align}\label{eq:DefFcts}
		\frac{\tilde\psi^\even_k}{\scriptstyle \sqrt{2}} := \begin{cases} 
     \sin(2 k \pi \cdot), \, k > 0, \\ 
    \cos(2 k \pi \cdot), \, k < 0,
\end{cases} 
		\hspace{1.4em}
		\frac{\tilde\psi^\odd_k}{\scriptstyle \sqrt{2}} := \begin{cases} 
    \sin( (2k-1)\pi \cdot), \, k > 0, \\ 
    \cos( (2k-1) \pi \cdot), \, k < 0.
\end{cases}
	\end{align}
Both, $\tilde\Psi^\even:=\{ \tilde\psi_k^\even :\, k\in \Z \}$ and $\tilde\Psi^\odd:=\{ \tilde\psi^\odd_k :\, k\in \mathbb{Z} \setminus \{0\} \}$ are orthonormal bases for $L_2(\Omega)$. The Fourier expansion of $h\in L_2(\Omega)$ then reads
\begin{align}\label{eq:FourierOmega}
	h = \sum_{k\in\Z} \hat{d}_k(h) e^{2k\pi i \cdot},
	\qquad\text{where} \quad
	\hat{d}_k(h) = \langle h, e^{2k\pi i \cdot}\rangle_{L_2(\Omega; \C)}
\end{align}
are the Fourier coefficients.\hfill$\diamond$
\end{compactenum}
\end{remark}
}

\rev{
We note the following simple fact for later reference.
\begin{lemma}\label{Lem:HWS}
	Let $h \in L_2^\even$, then $H(\tau) := \int_{\tau-1}^{\tau} h(s) ds=H(1)$ for all $\tau \in [0,1]$.
\end{lemma}
\begin{proof}
	Since  $h \in L_2^\even$, we get $\int_{\tau-1}^{0} h(s) ds = \int_{\tau}^{1} h(s-1) ds = \int_{\tau}^{1} h(s) ds$, so that we conclude that
$H(\tau) =  \int_{\tau-1}^{0} h(s) ds + \int_{0}^{\tau} h(s) ds  =  \int_{\tau}^{1} h(s) ds +  \int_{0}^{\tau} h(s) ds = H(1)$.
\end{proof}
}

\rev{
\begin{remark}\label{Rem:L2scaling}
We will also need to relate HWS functions on $\OP$ with their restriction to $\Omega$ making use of the half-wave symmetry. Let $g\in L_2^\even$, then it is easily verified that $\hat c_{2k}(g) = \sqrt{2}\, \hat d_{2k}(g)$, which means that $\| g\|_{L_2(\OP;\C)}^2 = 2 \| g\|_{L_2(\Omega;\C)}^2$. This latter relation also holds for  $g\in L_2^\odd$ and thus for all $g\in L_2^\HWS$.
\hfill$\diamond$
\end{remark}
}

\subsection{Sobolev spaces of broken order}\label{sobolev-spaces}
\rev{Our aim is to} relate the decay of the Kolmogorov $N$-width to the regularity of the parameter-dependent solutions. \rev{To this end, 
we define the $H^r(\OP)$-norm by the $H^r(\OP)$-seminorm\footnote{By first defining the fractional derivative $g^{(r)} \in L_2(\OP)$, the seminorm can also be defined as $\lvert g \rvert_{H^r(\OP)} := \| g^{(r)} \|_{L_2(\OP)}$. But we will not need the explicit (weak) fractional derivative for the estimates.} of broken order $r\geq0$ by 
\begin{align*}
\lvert g \rvert_{H^r(\OP)}^2 := \sum_{k \in \Z}^{\infty} (\pi k)^{2r} |\hat{c}_k(g)|^2  < \infty, \qquad \| g\|_{H^r(\OP)}^2 := \| g\|_{L_2(\OP)}^2 + \lvert g \rvert_{H^r(\OP)}^2,
\end{align*}
to define the Sobolev space of (broken) order by\footnote{It holds $H^{r} (\OP) 
	:= \Big\{ g \in L_2(\OP): \| g\|_{H^r(\OP)}^2 < \infty \Big\}= H^{r} (\OP) 
	:= \Big\{ g \in L_2(\OP): \lvert \, g \, \rvert_{H^r(\OP)}^2 < \infty \Big\}$.}
\begin{align}\label{eq:Hr}
	H^{r} (\OP) 
	&:= \Big\{ g \in L_2(\OP): \| g\|_{H^r(\OP)}^2 = \sum_{k \in \Z}^{\infty} (1+(\pi k)^{2r}) |\hat{c}_k(g)|^2  < \infty \Big\}.
\end{align}%
Next, we define the corresponding spaces for even and odd half-wave symmetric functions, i.e., $H^{r,\odd} := L_2^{\odd} \cap H^{r}(\OP)$, $H^{r,\even} := L_2^{\even} \cap H^{r}(\OP)$ and $H^{r,\HWS} := L_2^{\HWS} \cap H^{r}(\OP)$. Of course, there are also other definitions of $r$-th order Sobolev spaces in the literature, which are (often) equivalent to the above setting.} 

\rev{%
We shall also need the analogue on $\Omega$, i.e., 
\begin{align}\label{eq:HrOmega}
	H^{r} (\Omega) 
	&:= \Big\{ g \in L_2(\Omega): 
		\| g\|_{H^r(\Omega)}^2 = \sum_{k \in \Z}^{\infty} (1+(2\pi k)^{2r}) |\hat{d}_{k}(g)|^2  < \infty 
		\Big\}.
\end{align}%
\begin{corollary}\label{Cor:NormsHr}
	Let $g\in H^{r,\HWS}$, then 
	$$ \|g \|_{L_2(\OP)} = \sqrt{2}\, \|g \|_{L_2(\Omega)} , \quad \lvert g \rvert_{H^r(\OP)} = \sqrt{2}\, \lvert g \rvert_{H^r(\Omega)}, \quad \|g \|_{H^r(\OP)} = \sqrt{2}\, \|g \|_{H^r(\Omega)}.$$
\end{corollary}
\begin{proof}
	By Remark \ref{Rem:L2scaling}, we have $\hat c_{2k}(g) = \sqrt{2}\, \hat d_{2k}(g)$ and inserting this into \eqref{eq:Hr} and \eqref{eq:HrOmega} shows the claim.
\end{proof}
}

\section{\rev{Optimal shift-isometric spectral} decomposition}\label{section3}
Since the orthogonal projection is the best approximation, we are considering orthogonal decompositions of the spaces that are relevant for the transport problem. \rev{However}, the orthogonal space decompositions need to be tailored to a given function $g$ in order to bound or represent $d_N(\cU^g)$ and $\delta_N(\cU^g)$ for that function $g$. 

\subsection{Eigenfunctions and \rev{$L_2$-}optimality}\label{Sec:EigenSpaces}
\rev{It is not surprising that spaces} spanned by eigenfunctions are relevant for analyzing the $L_2$-average $N$-width $\delta_N(\cU^g)$. In fact, the appropriate orthogonal space decomposition is built by eigenspaces of the operator induced by the bilinear form
\begin{align*}
    k^g: L_2(\Omega)\times L_2(\Omega) \to \R,
    \qquad
    k^g(v,w) :=  \tfrac{1}{|\cP|} \int_{\cP} 
    \langle u_{\mu}, v \rangle_{L_2(\Omega)}\, 
    \langle u_{\mu}, w \rangle_{L_2(\Omega)}\, d \mu ,
\end{align*}
with \enquote{snapshots} $u_\mu=g(\cdot-\mu)$ induced by the function $g$ are defined by \eqref{eq:functiontobeapproximated}. We need to keep track on the dependence on $g$. Obviously, $k^g$ is a symmetric and positive semi-definite bilinear form. We define the induced operator by
\begin{align*}
    \cK^g: L_2(\Omega)\to L_2(\Omega),
    \qquad
    \langle \cK^g\varphi,\psi\rangle_{L_2(\Omega)} :=
    k^g(\varphi,\psi), \hspace{0.4em} \psi\in L_2(\Omega),
\end{align*}
which is consequently a positive semi-definite operator.

\begin{remark} \label{rem-kernel-expr}
    We note the following representation of $k^g$ (and hence $\cK^g$):
    \begin{align*}
        k^g(v,w)
        &= \int_\Omega \int_\Omega 
        v(x) 
        \kappa^g(x,y)
        w(y)\, dy\, dx,
    \end{align*}
    where $\kappa^g(x,y) := \langle u_\mu(x), u_\mu(y) \rangle_{L_2(\cP)} =  \tfrac{1}{|\cP|} \int_\cP g(x-\mu)\, g(y-\mu)\, d\mu$, i.e., $k^g$ is an integral operator with kernel $\kappa^g$.\hfill$\diamond$
\end{remark}

The operator $\cK^g$ admits an $L_2(\Omega)$-ON basis $\{ v^g_{ k}\}_{k\in\N}$ of eigenfunctions according to non-negative \textbf{ordered}\footnote{\rev{This will be relevant later.}} eigenvalues $\lambda^g_{1}\ge \lambda^g_{2}\geq\cdots\ge 0$.\footnote{$\cK^g v^g_{i}=\lambda^g_{i} v^g_{i}$, i.e., $k^g(v^g_{i},\psi)=\lambda^g_{i} \langle v^g_{i},\psi\rangle_{L_2(\Omega)}$, for all $i\in\N$.}
Then, we define
\begin{align} \label{defineVNbyeigs}
    V^g_N := \Span\{ v^g_1,...,v^g_N\}
\end{align}
along with the orthogonal projector $P_N^g: L_2(\Omega)\to V^g_N$ defined as $P_N^g v := \sum_{k=1}^N \langle v, v^g_k\rangle_{L_2(\Omega)} \, v^g_k$.  
By orthonormality \rev{and the definition of $k^g$, we get that sorted eigenfunctions are optimal.}

\rev{
\begin{lemma} \label{Eigs-are-optimal}
Let $N \in \N$. The linear space defined in \eqref{defineVNbyeigs} is optimal w.r.t.\ the $L_2$-average $N$-width, i.e.
\begin{align} \label{eq:L2opt}
	\delta_N(\cU^g)^2 
	&=\dist(V^g_N,\cU^g)_{L_2(\cP;L_2(\Omega))}^2  
	= \sum_{k=N+1}^{\infty} \lambda_k^g .
\end{align}
\end{lemma}
\begin{proof}
Let $H_{N} \subset L_{2}(\Omega)$ be some $N$-dimensional subspace generated by an ON basis $\{\varphi_{1},..., \varphi_{N}\}$. The Basis Extension Theorem allows us to extend to an ON basis $\{\varphi_{k}\}_{k \in \N}$ of $L_2(\Omega)$. We expand $\varphi_k$ in terms of the eigenbasis $\{ v_i^g\}_{i\in\N}$ as  $\varphi_k =  \sum_{i\in\N} \langle \varphi_k, v_i^g \rangle_{L_2(\Omega)} v_i^g,$ for all $k\in \N$. Then, by using the orthogonal projection onto $H_N$, 
\begin{align*}
	\dist(H_N,\mathcal{U}_g)_{L_2(\mathcal{P};L_2(\Omega))}^2 
	&=  \sum_{k=N+1}^{\infty} k^g(\varphi_k,\varphi_k)
	= \sum_{k=N+1}^{\infty} \sum_{i=1}^{\infty} 
	\langle \varphi_k, v_i^g \rangle_{L_2(\Omega)} \kern-5pt 
	\underbrace{k^g(v_i^g,\varphi_k)}_{=\lambda_i^g \langle v_i^g, \varphi_k\rangle_{L_2(\Omega)}}\\
	&=  \sum_{i=1}^{\infty}  \lambda_i^g  \sum_{k=N+1}^{\infty}  \langle v_i^g, \varphi_k\rangle_{L_2(\Omega)}^2.
\end{align*}
\vspace*{-12pt}
\begin{align*}
 \text{Next, } \quad \qquad	&\dist(H_N,\mathcal{U}_g)_{L_2(\mathcal{P};L_2(\Omega))}^2 
	- \dist(V_N^g,\mathcal{U}_g)_{L_2(\mathcal{P};L_2(\Omega))}^2 =  \\
	&\kern+20pt=  \tfrac{1}{|\cP|} \int_\cP \inf_{h_N\in H_N} \| u_\mu- h_N\|_{L_2(\Omega)}^2 d\mu
		-  \tfrac{1}{|\cP|} \int_\cP \inf_{w_N\in V_N^g} \| u_\mu-w_N\|_{L_2(\Omega)}^2 d\mu \qquad \qquad \qquad \\
	&\kern+20pt=\underbrace{  \sum_{i=1}^\infty \lambda_i^g 
		 \sum_{k=N+1}^\infty \langle\varphi_k, v_i^g\rangle_{L_2(\Omega)}^2}_{= \sum_{k=N+1}^\infty k^g(\varphi_k,\varphi_k)}
		- \underbrace{ \sum_{k=N+1}^\infty k^g( v_k^g , v_k^g )}_{=\sum_{k=N+1}^\infty \lambda_k^g} .
\end{align*}
By orthonormality, we have
$0 \leq \theta_i^N :=  \sum_{k=1}^{N} \langle \varphi_k, v_i^g \rangle_{L_2(\Omega)}^2 
	\leq  \sum_{k=1}^{\infty} \langle \varphi_k, v_i^g \rangle_{L_2(\Omega)}^2 
	= \| v_i^g \|_{L_2(\Omega)}^2 = 1$ and $N =  \sum_{k=1}^{N} \|\varphi_k\|_{L_2(\Omega)}^2 =  \sum_{k=1}^{N}  \sum_{i=1}^{\infty} \langle \varphi_k, v_i^g \rangle_{L_2(\Omega)}^2 =  \sum_{i=1}^{\infty} \theta_i^N$. Hence,
\begin{align*}
	 \sum_{k=N+1}^\infty k^g(\varphi_k,\varphi_k)
	&=  \sum_{i,k=1}^\infty \lambda_i^g \langle \varphi_k, v_i^g\rangle_{L_2(\Omega)}^2
	-  \sum_{k=1}^N  \sum_{i=1}^\infty \lambda_i^g \langle \varphi_k, v_i^g\rangle_{L_2(\Omega)}^2 \\
	&= \sum_{i=1}^\infty \lambda_i^g \underbrace{ \sum_{k=1}^\infty \langle \varphi_k, v_i^g\rangle_{L_2(\Omega)}^2}_{= \| v_i^g\|_{L_2(\Omega)}^2=1}
	-  \sum_{i=1}^\infty \lambda_i^g  \theta_i^N
	=  \sum_{i=1}^\infty (1-\theta_i^N) \lambda_i^g.
\end{align*}
Thus, finally, we get
\begin{align*}
	&\dist(H_N,\mathcal{U}_g)_{L_2(\mathcal{P};L_2(\Omega))}^2 
	- \dist(V_N^g,\mathcal{U}_g)_{L_2(\mathcal{P};L_2(\Omega))}^2 = \\
	&\kern+20pt=  \sum_{i=1}^\infty (1-\theta_i^N) \lambda_i^g
		-  \sum_{i=N+1}^\infty \lambda_i^g
	=  \sum_{i=1}^N (1-\theta_i^N) \lambda_i^g -  \sum_{i=N+1}^\infty \theta_i^N \lambda_i^g \\
	&\kern+20pt\ge \lambda_{N}^g \bigg(  \sum_{i=1}^N (1- \theta_i^N) 
		-  \sum_{i=N+1}^{\infty} \theta_i^N \bigg) 
		= \lambda_{N}^g \bigg( N -  \sum_{i=1}^{\infty} \theta_i^N \bigg)
		= 0,
\end{align*}
which proves the claim.
\end{proof}}

In this sense, eigenspaces are optimal in $L_2$, which is of course quite well-known from the Singular Value Decomposition (SVD) or the Proper Orthogonal Decomposition (POD) in MOR. 

\subsection{Shift-isometry \rev{and $L_{\infty}$-}optimality}\label{Sec:ShiftOthDec}
In order to link $\delta_N(\cU^g)$ to the $L_\infty$-based $N$-width $d_N(\cU^g)$, we need an additional property of an orthogonal decomposition to be introduced next. \rev{To this end, we shall need different bases for shift-isometric subspaces}.\footnote{In our case, all shift-isometric subspaces will have dimension 2 except for the space spanned by the constant function, which has dimension 1.}

\begin{definition} \label{ortho-subspace-basis}
	Let $g\in L_2(\OP)$. An orthogonal space decomposition of $L_2(\Omega)$ induced by a family of subspaces $\cW^g:=\{W_k^g\}_{k \in\N}$ with the associated orthogonal projectors $\cQ^g:=\{Q_k^g\}_{k \in\N}$, $Q_k^g: L_2(\Omega)\to W_k^g$, is called \emph{shift-isometric orthogonal decomposition} (w.r.t.\ $g$) of $L_2(\Omega)$ if               \begin{align}\label{eq:shiftinvariant}
                  \big\| Q_k^g g(\cdot-\mu) \big\|_{L_2(\Omega)} = \big\| Q_k^g g \big\|_{L_2(\Omega)}
                  \quad\text{for all } \mu\in\cP ,
              \end{align}
        i.e., if the orthogonal projectors are \emph{shift-isometric}. \hfill$\diamond$
\end{definition}%

\begin{remark}
    Note, that $g$ needs to be defined on the larger space $\OP$ in order to apply $g(\cdot-\mu)$.
    However, the solution of the transport problem \eqref{eq:transport} is defined on the domain $\Omega$.
    Whenever we take a norm $\|\cdot\|_{L_2(\Omega)}$ or apply $Q_k^g$, we consider implicitly the restrictions $g_{|\Omega}$ or $g(\cdot- \mu)_{|\Omega}$ respectively.\hfill$\diamond$
\end{remark}

Now, we start by \emph{assuming} that such a \rev{shift-isometric orthogonal decomposition} $\cW^g$ exists and show that \eqref{eq:shiftinvariant} is a key property. If $\cW^g$ is a (shift-isometric) orthogonal decomposition of $L_2(\Omega)$, each $u_\mu\in\cU^g\subset L_2(\Omega)$ has a unique decomposition $u_\mu = \sum_{k=1}^{\infty} Q_{k}^g u_\mu$. Let \rev{$M \in \N$}, we define
\begin{align}\label{eq:XN}
    \rev{X}_N^g := \bigoplus_{k=1}^{M} W_k^g,
    \hspace{1.3em} \text{where} \hspace{0.6em}
    N := N(M)=\dim (\rev{X}_N^g) = \sum_{k=1}^{M} \dim(W_k^g)
\end{align}
as a candidate for the best approximation space in the sense of Kolmogorov.\footnote{Note, that by fixing the dimensions of $W_k^g$ a priori, one might not always be able to construct spaces $V_N^g$ of any dimension $N\in\N$, since $N$ must be a sum of the dimensions of the $W_k^g$, $k=1,...,M(N)$. As an example, if $\dim(W_k^g)=2$ for all $k \in \mathbb{N}$, $N$ must be even.} Later, the spaces $W^g_k$ will be spanned by two eigenfunctions corresponding to the same eigenvalue.
Clearly, the approximation
$\tilde v_N := \sum_{k=1}^{M(N)} Q_k^g u_\mu$ converges to $u_\mu$ as $N\to\infty$, which implies that both $d_N( \cU^g )$ and $\delta_N( \cU^g)$ converge towards zero. Moreover, the orthogonality easily allows us to control the error.

\begin{proposition}\label{Prop:distProj}
    Let $g\in L_2(\OP)$, let $\cW^g$ be a corresponding shift-isometric orthogonal decomposition of $L_2(\Omega)$  and let $\rev{X}_N^g$ be defined as in \eqref{eq:XN}. Then, 
    \begin{align}
        \dist(\rev{X}_N^g,\cU^g)_{L_\infty(\cP;L_2(\Omega))}^2
         & \kern-2pt=   \dist(\rev{X}_N^g,\cU^g)_{L_2(\cP;L_2(\Omega))}^2
        \kern-2pt=\kern-10pt  \sum_{k=M(N)+1}^{\infty}  \kern-10pt\| Q_{k}^g g \|_{L_2(\Omega)}^2.
        \label{eq:dist}
    \end{align}
\begin{proof}
    The second equality originates from the orthogonal decomposition. Concering $L_\infty(\cP)$, we get by shift-isometry (Def.~\ref{ortho-subspace-basis}) 
    \begin{align*}
    	\| Q_k^g u_\mu\|_{L_2(\Omega)}=\| Q_k^g g(\cdot-\mu)\|_{L_2(\Omega)}=\| Q_k^g g\|_{L_2(\Omega)}, 
	\end{align*}
	\rev{(and this is the key property to eliminate the $\mu$-dependence)} so that 
    \begin{align*}
        \dist(\rev{X}_N^g, \cU^g)_{L_\infty(\cP;L_2(\Omega))}^2
          = \sup_{\mu\in\cP} \| u_\mu - P_N^g u_\mu\|_{L_2(\Omega)}^2 = \sup_{\mu\in\cP} \sum_{k=M(N)+1}^{\infty} \| Q_k^g g\|_{L_2(\Omega)}^2         \\
         = \sum_{k=M(N)+1}^{\infty} \kern-10pt\| Q_k^g g\|_{L_2(\Omega)}^2 = \tfrac{1}{|\cP|} \int_\cP \sum_{k=M(N)+1}^{\infty} \| Q_k^g g\|_{L_2(\Omega)}^2 d\mu
        = \dist(\rev{X}_N^g, \cU^g )_{L_2(\cP;L_2(\Omega))}^2,
    \end{align*}
    which completes the proof.
\end{proof}
\end{proposition}

\rev{Recall, that the eigenvalues $\lambda_k^g$ are assumed to be sorted in decreasing order, which will be important below. In order to characterize the error of an approximation, we also need to keep track of the multiplicities of multiple eigenvalues. To this end, we introduce an enumeration by $\rho : \N \to \N$ so that\footnote{For the corollary, this can be weakened to $\lambda_{\rho(j)}^g \geq \lambda_{\rho(j+1)}^g ,$ $j \in \N$.}}
\begin{align*}
	\lambda_{\rho(1)}^g > \lambda_{\rho(2)}^g 
	> \lambda_{\rho(3)}^g 
	> \cdots 
	\quad \text{with} \quad 
	\lambda_{\rho(k)}^g = \lambda_{\rho(k)+1}^g = \dots = \lambda_{\rho(k+1)-1}^g, k\in \N,
\end{align*} 
and $m_k= \rho(k+1)-\rho(k)$ being the algebraic multiplicity of $\lambda^g_{\rho(k)}$, $k \in \N$. We summarize our findings.

\rev{
\begin{corollary} \label{cor-Nwidths-equal}
Let $M \in \N$ and $N=N(M) = \sum_{k=1}^{M} m_k \in \N$. Define $V_N^g := \bigoplus_{k=1}^{M} W_{\rho(k)}^g$, where $W_{\rho(k)}^g := \operatorname{span} \big\{ v_{\rho(k)}^g, v_{\rho(k)+1}^g, \dots,v_{\rho(k+1)-1}^g \big\}$ is the eigenspace corresponding to the eigenvalue $\lambda_{\rho(k)}^g$. 
If $\mathcal{W}^g= \big\{ W_{\rho(k)}^g \big\}_{k \in \N }$ is a shift-isometric decomposition, then 
\vspace*{-12pt}
\begin{align*}
	\delta_N(\cU^g)=d_N(\cU^g) = \Big(  \sum_{k=M+1}^{\infty} m_k \, \lambda^g_{\rho(k)} \Big)^{1/2} . 
\end{align*}
\end{corollary}
\begin{proof}
From Lemma \ref{Eigs-are-optimal} and Proposition \ref{Prop:distProj}, we have $d_N(\cU^g) \leq \dist(V_N^g, \cU^g)_{L_{\infty}(\cP;L_2(\Omega))} = \dist(V_N^g, \cU^g)_{L_2(\cP;L_2(\Omega))} = \delta_N( \cU^g ) \leq d_N( \cU^g)$.
\end{proof}
}

\section{\textit{N}-widths for half-wave symmetric functions} \label{Sec:NWidth}
We will now specify on the linear transport problem with the setting from section \ref{Subsec:Fourier} with $I=(0,1)$, $\Omega=(0,1)$, $\cP=[0,1]$ and $\OP=(-1,1)$. We are going to construct shift-isometric \rev{spectral} decompositions for half-wave symmetric functions $g\in L_2^\HWS$ in terms of trigonometric functions. \rev{Plus,} these spaces $W_k^\HWS\subset L_2(\Omega)$ will have an additional property, namely they are \emph{shift-invariant}, i.e., \rev{$w \in \operatorname{span} \big\{ \xi_1, \xi_2 \big\} \subset L_2(\OP)$ implies that $w(\cdot - \mu)|_{\Omega} \in \operatorname{span} \big\{ \xi_1|_{\Omega}, \xi_2|_{\Omega} \big\}$ for all $\mu \in \mathcal{P}$}. We will need to consider even and odd HWS functions separately.

\subsection{\rev{Trigonometric shift-isometric spectral decompositions}}\label{sect:eig-decomp}
\smallskip
\subsubsection*{\rev{Eigenspaces}}
\rev{Let us start by reporting some properties of the basis in \eqref{eq:DefFcts} for even HWS functions. The odd ones will be considered afterwards.}

\begin{lemma}\label{La:Weven}
     For any $k \in \mathbb{N}$, the set \rev{$\{\tilde\psi_k^\even, \tilde\psi_{-k}^\even\} \subset {L_2(\Omega)}$}\footnote{See Remark \ref{Rem:L2Omega}.} is an orthonormal basis for $\Wke:=  \Span \{\tilde\psi_k^\even, \tilde\psi_{-k}^\even\} = \rev{ \{ f:\Omega\to\R: f(x) = c\, e^{2k\pi ix} + \bar{c}\, e^{-2k\pi ix}, c\in\C\} }$ with $\dim(\Wke) = 2$ and these spaces are shift-invariant.
\begin{proof}
    The statements concerning orthonormality and dimension are straightforward.
    Let $\mu\in\cP$, then for $k\in \N$
    \rev{$\tilde\psi_{-k}^\even(\cdot - \mu)|_{\Omega}
        = \cos ( 2k\pi \mu ) \, \sqrt{2} \cos ( 2k\pi \cdot )
            + \sin ( 2k\pi \mu ) \, \sqrt{2} \sin ( 2k\pi \cdot ) \in \Wke$,
    and the} same applies for $\tilde\psi_k^\even$.
\end{proof}
\end{lemma}
\rev{A} simple proof shows that  $\Wke = \{ \alpha \,\sin(2k\pi \cdot + \beta):\, \alpha, \beta \in \R \}$.
For $k=0$, we set $W_0^{\mathrm{\even}}:= \operatorname{span} \{ 1 \}, \dim(W_0^{\mathrm{\even}}) = 1$, i.e., the constant functions.

\begin{remark}
    The above definition is similar to Kolmogorov's paper in 1936 \cite{kolmogoroffUberBesteAnnaherung1936}, where the best basis functions for all periodic $H^r(0,1)$-functions, $r \in \N_0$, with $\|f^{(r)}\|_{L_2(0,1)} \leq 1$ is shown to be $\{1$, $\sqrt{2 } \sin \left( 2 \pi k \, \cdot \, \right)$, $\sqrt{2}  \cos \left( 2 \pi k \, \cdot \, \right)$,  $k = 1, 2,..., \tfrac{N-1}{2}\}$. For such classes of functions, Kolmogorov quantified a constant, which was later called \emph{Kolmogorov $N$-width} in honor of his contributions and \rev{he also} proved $d_{N} = (\pi N)^{-r}$.\hfill$\diamond$
\end{remark}

\rev{
\begin{lemma}\label{Lem:EVeven}
 Let $g\in L_2^\even$. Then, $\tilde\psi_{\pm k}^\even$, $k\in\N_0$, are $L_2(\Omega)$-normalized eigenfunctions of $\cK^g$  corresponding to the eigenvalue \footnote{Recall that  $\hat{c}_{-2k}(g)=\overline{\hat{c}_{2k}(g)}$, so that their absolute values coincide.}
	\begin{align}\label{eq:EVeven}
		\tilde\lambda_{\pm k}^{g}
		:= \tfrac12 \, |\hat{c}_{2k}(g)|^2
		= \tfrac12 \, |\hat{c}_{-2k}(g)|^2.
	\end{align}
	\end{lemma}
\begin{proof}
	    The last equality follows from \eqref{eq:ckg}. Next, we insert \eqref{eq:Fourier-hws} into the kernel $\kappa^g$ from Remark \ref{rem-kernel-expr} and get for $g\in L_2^\even$
	\begin{align*}
		\kappa^g(x,y) 
		&= \int_{\mathcal{P}} g(x-\mu) g(y-\mu) d \mu = \int_{\cP} \Big[ \sum_{k \in \Z} \hat{c}_{2k}(g) \, 
			\tfrac{\textstyle e^{2i \pi k (x-\mu)}}{\sqrt{2}} 
			\sum_{\ell \in \Z} \hat{c}_{2\ell}(g) \, 
			\tfrac{\textstyle e^{2i \pi \ell (y-\mu)}}{\sqrt{2}} 
			\Big] d \mu \\
	 &= \tfrac12 \sum_{k,\ell \in \Z} \hat{c}_{2k}(g)
	 	\hat{c}_{2\ell}(g) e^{2 \pi i (kx + \ell y) } \delta_{k,-\ell} 
	= \tfrac12 \sum_{\ell \in \Z} |\hat{c}_{2\ell}(g)|^2 e^{2 \pi i \ell (x-y) }. 
	\end{align*}
	Hence, for any $w \in L_2(\Omega)$, 
	\begin{align*}
	k^g( \tilde\psi_{- k}^\even , w ) 
	&= \int_\Omega \int_\Omega \tilde\psi_{- k}^\even(x) \, \kappa^g(x,y)\, 
		w(y)\, dx\, dy \\
	& = \tfrac12 \sum_{\ell\in\Z} |\hat c_{2\ell}(g)|^2 \int_\Omega\int_\Omega 
		 \tilde\psi_{- k}^\even(x)\, e^{2\pi i\ell(x-y)} w(y)\, dx\, dy \\
	&= \tfrac12 \sum_{\ell\in\Z} |\hat c_{2\ell}(g)|^2 
		\int_\Omega \underbrace{\int_\Omega \sqrt{2}\, \cos(2k\pi x)\, e^{2\pi i \ell x} dx}_{=2^{-1/2}(\delta_{k,\ell}+\delta_{k,-\ell})} e^{-2\pi i\ell y}\, w(y)  dy \\
	&= \tfrac1{2\sqrt{2}} |\hat c_{2k}(g)|^2 \int_\Omega w(y) e^{-2\pi i ky} dy
		+ \tfrac1{2\sqrt{2}} |\hat c_{-2k}(g)|^2 \int_\Omega w(y) e^{2\pi i ky} dy \\
	&=  \tfrac1{2\sqrt{2}}  |\hat c_{2k}(g)|^2
	\int_\Omega w(y) ( e^{2\pi i ky}+e^{-2\pi i ky})\, dy \\
	&= \tfrac12 \, |\hat{c}_{2k}(g)|^2 
		\langle \tilde\psi_{- k}^\even , w \rangle_{L_2(\Omega)} = \tilde\lambda_{-k}^{g} \langle \tilde\psi_{- k}^\even , w \rangle_{L_2(\Omega)},
	    \end{align*}
	    which proves the claim for $\tilde\psi_{- k}^\even$. A similar derivation applies to $\tilde\psi_{k}^\even$.
\end{proof}
}
In a quite analogous manner, we get \rev{similar results} for the odd \rev{HWS}  case. We skip the proofs.

\begin{lemma}\label{La:Wodd}
     For any $k \in \mathbb{N}$, the set $\{ \rev{\tilde\psi_{\pm k}^\odd}\}\subset {L_2(\Omega)}$ is an orthonormal basis of the shift-invariant spaces $\Wko:=  \Span \{ \rev{\tilde\psi_{\pm k}^\odd} \}$ with $\dim(\Wko) = 2$. \hfill\qed
\end{lemma}

\rev{
\begin{lemma}\label{La:EVodd}
Let $g\in L_2^\odd$. The functions $\rev{\tilde\psi_{\pm k}^\odd}$, $k\in\N$,  are $L_2(\Omega)$-normalized eigenfunctions of $\cK^g$ corresponding to the eigenvalue 
    \begin{align}\label{eq:EVodd}
        \tilde\lambda_{\pm k}^{g} 
        := \tfrac12 \, |\hat{c}_{2k-1}(g)|^2 
        = \tfrac12 \, |\hat{c}_{-(2k-1)}(g)|^2 .
    \end{align}
    \hfill\qed
\end{lemma}
    }
 \vspace*{-0.2cm}

\smallskip
\subsubsection*{\rev{Shift-isometry}}
We shall now prove that the above construction yields shift-isometric orthogonal decompositions. It turns out that shift-isometry (Def.~\ref{ortho-subspace-basis}) and half-wave symmetry allow us to use the same orthogonal decomposition of $L_2(\Omega)$ for all $L_2^\even$- and $L_2^\odd$-functions, respectively. We do not need a specific decomposition $\cW^g$ for each individual function $g \in L_2^\HWS$.

\begin{lemma}\label{Lemma-ehws-shiftiso}
    \begin{compactenum}[(a)]
        \item The family $\cW^\even:=\{\Wke \}_{k \in\N_0}$ is a shift-isometric orthogonal decomposition of $L_2(\Omega)$ w.r.t.\ all $g\in L_2^\even$.
        \item The family $\cW^\odd:=\{\Wko \}_{k \in\N}$ is a shift-isometric orthogonal decomposition of $L_2(\Omega)$ w.r.t.\ all $g\in L_2^\odd$.
    \end{compactenum}
\begin{proof}
	We restrict ourselves to the even case \textit{(a)} and note that \textit{(b)} is analogous. Recall the Fourier expansion \eqref{eq:Fourier-hws} of $g \in L_2^\even$, which yields
	\begin{align*}
	    g (\cdot - \mu)                 
	    &=  \sum_{k \in \Z}  \hat{c}_{2k}(g ) 
	    \tfrac{\textstyle e^{2 i \pi k ( \cdot - \mu )}}{\sqrt{2}},
	    \qquad\text{so that} \\
	    [Q_k^\even g (\cdot- \mu)](\bullet)
	    &= \langle g (\cdot- \mu),   e^{2 \pi i k \cdot} \rangle_{L_2(\Omega,\C)}\,
	    		 e^{2 \pi i k \bullet} 
	    	+ \langle g (\cdot- \mu),  e^{-2 \pi i k \cdot} \rangle_{L_2(\Omega,\C)}\, 
			e^{-2 \pi i k \bullet} \\
	    &= \tfrac{ \hat{c}_{2k}(g )}{\sqrt{2}} e^{2 \pi i k ( \bullet -\mu ) } +\tfrac{ \hat{c}_{-2k}(g )}{\sqrt{2}} e^{-2 \pi i k ( \bullet -\mu ) } = [Q_k^\even g](\bullet-\mu).
	    \end{align*}
	    Now, we apply Lemma \ref{Lem:HWS} to $h(\bullet) :=\lvert [Q_k^\even g] (\bullet) \rvert^2$, noting that $h(x+1)=\lvert [Q_k^\even g](x+1) \rvert^2= \lvert[Q_k^\even g](x) \rvert^2=h(x)$ for $x\in [-1,0]$ a.e., i.e., $h\in L_2^{\even}$. Then, by Lemma \ref{Lem:HWS}
	    \[
\| Q_k^\even g (\cdot - \mu)\|_{L_2(\Omega)}^2
= \| [Q_k^\even g] (\bullet - \mu)\|_{L_2(\Omega)}^2
= \int_0^1 \kern-5pt h(s - \mu) \, ds
= \int_0^1 \kern-5pt h(s) \, ds
= \| Q_k^\even g \|_{L_2(\Omega)}^2
\]
	    for all $\mu\in\cP$, which proves the claim.
\end{proof}
 \end{lemma}
\subsection{Optimal sorting and size of the $N$-widths}
\rev{
Now, we are going to construct optimal subspaces $V_N^g$ in the sense of Kolmogorov for a given $g\in L_2^\HWS$ in terms of the above eigenspace decompositions. To this end, we have to deal with (at least) two technical difficulties, namely
\begin{compactitem}
	\item the eigenvalues determined in Lemma \ref{Lem:EVeven} and \ref{La:EVodd}  are related to the Fourier coefficients of $g$, which are in general not sorted;
	\item for even HWS functions, we have $\dim(W_0^{\even})=1$, whereas all other spaces are of dimension $2$; hence by collecting these spaces, we obtain $V_N^g$ with $N$ being odd for $g\in L_2^\even$ and $N$ even for $g\in L_2^\odd$. 
\end{compactitem} 
Let us start by introducing an appropriate sorting, where, due to the mentioned technicalities, we would need to distinguish between the even and odd HWS-case. To this end, we shall position the constant component separately later. Hence, denoting the mean of $g$ by $\bar{g} := \int_{\OP} g(x)\, dx = 0$, we have that $\bar{g}=0$ for $g\in L_2^\odd$. We set
\begin{align*}
	L_{2,0}^\even := \{ g\in L_2^\even:\, \bar{g}=0\} = L_2^\even/\R
\end{align*}
and --for a unified notation only-- $L_{2,0}^{\odd} := L_{2}^{\odd}$  as well as $L_{2,0}^\HWS := L_{2,0}^{\even} \cup L_{2,0}^{\odd}$.}

\rev{%
As for the first bullet above, recall from \S\ref{Sec:EigenSpaces} that the eigenvalues $\lambda_k^g$ of $\cK^g$ need to be sorted. By Lemma \ref{Lem:EVeven} and \ref{La:EVodd} we know that the eigenvalues $\tilde\lambda_k^g$ are related to the Fourier coefficients of $g$, but they are not sorted (which is also the reason for the \enquote{$\tilde{\,\,\,\,}$}-notation). Hence, we need an optimal sorting $\sigma: \N \to \N$ in the sense that 
\begin{align} \label{opt-sorting}
	\tilde\lambda^{g}_{\sigma(k)} 
	= \tilde\lambda^{g}_{-\sigma(k)} 
	\geq \tilde\lambda^{g}_{\sigma(k+1)} 
	= \tilde\lambda^{g}_{-\sigma(k+1)}, 
	\qquad
	k \in \N.
\end{align}
The first eigenvalue $\tilde\lambda_0^g$ corresponding to the constant function will be considered separately. Then, we set eigenvalue-eigenfunction pairs
$$
	\big( \lambda^{g}_{2k-1}, v^{g}_{2k-1} \big) := \big( \tilde\lambda^{g}_{\sigma(k)}, \tilde\psi^{\HWS}_{\sigma(k)} \big) , \quad \quad \big( \lambda^{g}_{2k}, v^{g}_{2k} \big) := \big( \tilde\lambda^{g}_{-\sigma(k)}, \tilde\psi^{\HWS}_{-\sigma(k)} \big) ,
	$$ 
which gives $\lambda_1^g = \lambda_2^g \ge \lambda_3^g = \lambda_4^g\ge \lambda_5^g \, \cdots$ as well as 
\begin{align}\label{eq:sortedlambda}
	 \lambda_{2k-1}^g = \lambda_{2k}^g 
	 = \begin{cases} 
		\tfrac{1}{2} \left| \hat{c}_{2\sigma(k)}(g) \right|^2 , & \text{if } g \in L_{2,0}^\even, \\[5pt]
		\tfrac{1}{2} \left| \hat{c}_{2\sigma(k)-1}(g) \right|^2 , & \text{if } g \in L_{2,0}^\odd.
	\end{cases} 
\end{align}
Recalling that each of the spaces $W_{k}^\HWS = \Span \{ \tilde\psi^{\HWS}_{k} , \tilde\psi^{\HWS}_{-k} \},$ $k \in \N,$ has dimension $2$, we set for $M\in\N$
\begin{align} \label{VN-def}
	V_{2M}^g := \bigoplus_{k=1}^{M} W_{\sigma(k)}^\HWS = \Span \{ v^{g}_{1}, v^{g}_{2}, ..., v^{g}_{2M-1}, v^{g}_{2M} \},
	\quad\text{i.e.,}\quad
	\dim(V_{2M}^g)=2M.
\end{align}
}
\rev{%
\begin{theorem}\label{Theorem:Optimal-simple}
Let $g\in L_{2,0}^{\HWS}$.
\begin{compactenum}[(a)]
	\item If $N\in\N$ is even, then $V_N^g$ is an optimal space of dimension $N$ in the sense of Kolmogorov with\ $d_N( \cU^g )^2 = \delta_N( \cU^g )^2 = \sum_{\ell = N+1}^{\infty} \lambda_{\ell}^{g}$. 
	\item If $N\in\N$ is odd, then $V_N^g:=V_{N-1}^g\oplus\Span\big\{\tilde\psi^\HWS_{\lceil \sigma(N/2) \rceil } \big\}$ and $V_N^g:=V_{N-1}^g\oplus\Span\big\{\tilde\psi^\HWS_{-\lceil \sigma(N/2) \rceil } \big\}$ satisfy $\dist(V_N^g, \cU^g )_{L_2(\cP;L_2(\Omega))} = \delta_N( \cU^g )$ as well as $\dist(V_N^g, \cU^g )_{L_\infty(\cP;L_2(\Omega))} = d_{N-1}( \cU^g )$, i.e., both spaces are $L_2(\mathcal{P})$-optimal, but do not improve the $L_{\infty}(\mathcal{P})$-distance.
	\end{compactenum}
\end{theorem}
\begin{proof}
	(a) The functions $\tilde\psi_{\pm k}^\HWS$ are eigenfunctions of $\cK^g$ by Lemma \ref{Lem:EVeven} and \ref{La:EVodd}. Moreover, $\{\lambda^g_{2k},\lambda^g_{2k-1}\}_{k\in\N}$  is the sequence of sorted eigenvalues corresponding to $\tilde\psi^\HWS_{\pm \sigma(k)}$. By Lemma \ref{Lemma-ehws-shiftiso}, $\cW^\HWS$ is a shift-isometric orthogonal decomposition of $L_2(\Omega)$ w.r.t.\ $g\in L_{2,0}^\HWS$, so that (a) follows from  Corollary \ref{cor-Nwidths-equal}.\\
	\smallskip
	To show \textit{(b)}, let $H_N = V_{N-1}^g\oplus\Span \{ \varphi \}$, where $\varphi \perp V_{N-1}^g$ and $\| \varphi\|_{L_2(\Omega)}=1$. Then, denoting by $P_{N-1}: L_2(\Omega)\to V_{N-1}^g$ the orthogonal projection onto  $V_{N-1}^g$, we have
	\begin{align*}
	\dist(H_N, \cU^g )_{L_\infty(\cP;L_2(\Omega))}^2 
	&= \sup_{\mu \in \cP} \big\| u_\mu - P_{N-1} u_{\mu} - \langle u_{\mu}, \varphi \rangle_{L_2(\Omega)} \varphi\big\|_{L_2(\Omega)}^2\\
	&\kern-100pt= \sup_{\mu \in \cP} \Big\{ \|u_{\mu}\|_{L_2(\Omega)}^2 - \|P_{N-1} u_{\mu}\|_{L_2(\Omega)}^2  - \langle u_{\mu}, \varphi \rangle_{L_2(\Omega)}^2  \Big\}  \\
	&\kern-100pt= \|u_{0}\|_{L_2(\Omega)}^2 - \|P_{N-1} u_{0}\|_{L_2(\Omega)}^2 
		- \inf_{\mu \in \cP} \langle u_{\mu}, \varphi \rangle_{L_2(\Omega)}^2 ,
	\end{align*}
	since $\cW^\HWS$ is shift-isometric and $ \|u_{0}\|_{L_2(\Omega)}^2 =  \|u_{\mu}\|_{L_2(\Omega)}^2$ by Lemma \ref{Lem:HWS} and $g^2 \in L_2^{\even}$. \\
	Next, we show that $\inf_{\mu \in \cP} \langle u_{\mu}, \varphi \rangle_{L_2(\Omega)}^2 = 0$. For $g \in L^{\even}_{2,0}$ (the case $g \in L^{\odd}_{2,0}$ is completely analogous), it holds that
	 $ 2 \int_0^1  g(x-\mu ) d \mu = \int_{x-1}^{x+1}  g(x- s ) ds = \int_{-1}^{1}  g(y) d y  =0$
	 for every $x \in (0,1)$ and therefore
	$\int_0^1 \int_0^1  g(x-\mu ) \varphi(x) d x \, d \mu = \int_0^1 \varphi(x) \int_0^1  g(x-\mu ) d \mu \, d x = 0$. 
	 Hence, by the mean value theorem, there exists $\bar \mu \in [0,1]$ such that 
	 $\int_0^1  g(x- \bar\mu ) \varphi(x) d x = 0$. 
	 \begin{align*}
	 \text{Hence,} \qquad 0 \le \inf_{\mu \in \cP} \langle u_{\mu}, \varphi \rangle_{L_2(\Omega)}^2
		&\le \langle u_{\bar\mu}, \varphi\rangle_{L_2(\Omega)}^2
		= \left(\int_0^1  g(x- \bar\mu ) \varphi(x) dx\right)^2 
		= 0. \qquad \qquad \qquad 
	 \end{align*}
	 As a result, doing the above steps backwards 
	 \begin{align*}
		\dist(H_N, \cU^g )_{L_\infty(\cP;L_2(\Omega))}^2 
		&= \|u_{0}\|_{L_2(\Omega)}^2 - \|P_{N-1} u_{0}\|_{L_2(\Omega)}^2\\
		&\kern-50pt= \sup_{\mu \in \cP} \Big\{ \|u_{\mu}\|_{L_2(\Omega)}^2 - \|P_{N-1} u_{\mu}\|_{L_2(\Omega)}^2   \Big\}
		= \dist(V_N^g, \cU^g )_{L_\infty(\cP;L_2(\Omega))}^2. 
	 \end{align*}
\end{proof}
}
\rev{
\begin{remark} \label{rem:dNdN-1}
As an immediate consequence of the $N$-widths definiton, it always holds $\delta_N( \cU^g ) \leq d_N( \cU^g ) \leq d_{N-1}( \cU^g ),$ $N \in \N$. Regarding Theorem \ref{Theorem:Optimal-simple} (a) for $g\in L_{2,0}^{\HWS}$ and even $N\in\N$, we have $\delta_N( \cU^g ) = d_N( \cU^g )$. \\
Regarding Theorem \ref{Theorem:Optimal-simple} (b) for odd $N\in\N$, we have
$$
	\delta_N( \cU^g )^2 
	= \sum_{\ell= N+1}^{\infty}\lambda_{\ell}^{g} \leq d_{N}( \cU^g )^2 \leq \dist(V_N^g, \cU^g )_{L_\infty(\cP;L_2(\Omega))} = d_{N-1}( \cU^g )^2 
	=  \sum_{\ell=N}^{\infty} \lambda_{\ell}^{g},
	$$
	which motivates to conjecture $d_{N}( \cU^g )= d_{N-1}( \cU^g )$ but it not proven, since in general terms as we cannot assume that any shift-isometric space is spanned by 2 basis functions. In fact, we found an example of a space spanned by $3$ linearly independent functions (which are not trigonometric functions) for which numerical computations indicate shift-isometry. \hfill$\diamond$
 \end{remark}
}
\medskip
\begin{example}[Discontinuous jump] \label{exam-jump-L2}
    We detail the example already investigated in \cite{ohlbergerReducedBasisMethods2016c}, where $d_N( \cU^g )\ge \tfrac{1}{2} N^{-1/2}$ was shown for discontinuous initial and boundary conditions, i.e., $g=\operatorname{sgn}(x)$, which is easily seen to be odd HWS.
    
Since $\langle g , \cos( (2k-1) \pi \cdot) \rangle_{L_2(\OP)} = 0$ and $\langle g , \sin((2k-1) \pi \cdot) \rangle_{L_2(\OP)} = \frac{\rev{4}}{(2k-1)\pi}$, we get $\rev{\lambda_{2k}^g =\lambda_{2k-1}^g } =\frac{4}{\pi^2(2k-1)^2}$ and $\lambda_k^g \ge \lambda^g_{k+1}$ for all $k\in\N$, i.e., sorting is not needed. Theorem \ref{Theorem:Optimal-simple} yields an exact representation of the $N$-width by
    \vspace*{-0.3cm}
        \begin{align}
        \delta_N( \cU^g )^2
            &= \sum_{k=N+1}^{\infty} \lambda_{2k-1}^g
            =\frac{4}{\pi^2}
                \sum_{k=N+1}^{\infty}
                \Big( 2\Big\lfloor \tfrac{k+1}{2} \Big\rfloor -1\Big)^{-2}
                    \nonumber \\
            &=  \frac{1}{\pi^2} \Psi^{(1)} 
                \Big( \left\lfloor \tfrac{N}{2} \right\rfloor 
                    + \tfrac{1}{2} \Big) 
                + \frac{1}{\pi^2} \Psi^{(1)} 
                    \Big( \left\lfloor \tfrac{N+1}{2} \right\rfloor 
                    + \tfrac{1}{2}  \Big) ,
        \label{equ:exact_decay_HS}
    \end{align} 
    where $\Psi^{(1)}(\cdot)$ is the first derivative of the Digamma function $\Psi(\cdot)$. Moreover, for even $N$, we have $\delta_N( \cU^g)=d_N( \cU^g )$.
    \demodrei
\end{example}
\rev{
\subsubsection*{Functions with a non-zero mean.}
So far, we assumed that $g\in L_{2,0}^\HWS$, i.e., $\bar g=0$. As $L_{2,0}^\odd=L_2^\odd$, this is no restriction for the odd HWS case. Hence, let us now consider the case $g\in L_{2}^{\even}$  with $\bar g\ne 0$ and hence $\hat c_0(g) = \tfrac1{\sqrt{2}} \bar g \ne 0$ as well as $\tilde\lambda^g_0 = \tfrac12 |\hat c_0(g)|^2=\tfrac14 \bar g^2\ne 0$.\\
Let $\lambda^g_1\ge\cdots\ge\lambda^g_N$ be the sorted $N$ largest eigenvalues of $\cK^g$. We have to distinguish two cases, namely
\begin{compactenum}[{case} 1:]
	\item $\lambda^g_N> \tilde\lambda^g_0$, i.e., the constant part is the smallest one, and
	\item $\lambda^g_N\le \tilde\lambda^g_0$, i.e., the constant part is significant in the sense that the eigenvalue corresponding to the constants $\tilde\lambda_0^{g}$ is \enquote{somewhere} in the ordering of the largest eigenvalues. 
\end{compactenum}
Then, the following table indicates how an optimal subspace can be chosen.
\begin{center}
\begin{tabular}{|c||c|c|c|}\hline
	&\text{} & \multicolumn{2}{c|}{\text{optimal space w.r.t.}} \\ 
	$N$ & case  & $L_2(\mathcal{P})$ / $\delta_N(\cU^g)$ & $L_{\infty}(\mathcal{P})$  / $d_N(\cU^g)$\\ \hline \hline
	\text{even} & 1 & $V_N^g$ & $V_N^g$ \\ \hline
	\text{even} & 2 & $V_{N-2}^g\oplus \Span \{ 1 \} \oplus\Span\{ \tilde\psi_{\lceil \sigma(N/2) \rceil }^\even \}$ & $\#$ \\ \hline
	\text{odd} & 1 & $V_{N-1}^g \oplus\Span\{ \tilde\psi_{\lceil \sigma(N/2) \rceil  }^\even \}$ & $\#$ \\ \hline
	\text{odd} & 2 & $V_{N-1}^g\oplus \Span \{ 1 \}$ & $V_{N-1}^g\oplus \Span \{ 1 \}$ \\ \hline
\end{tabular}
\end{center}\vspace*{5pt}
In the above table, \enquote{$\#$} in the last column means that those spaces are not clear, since as in Remark \ref{rem:dNdN-1}, we don't know which spaces are optimal when we don't have the shift-isometry property for the space spanned by the sorted eigenfunctions. For the other cases we know the optimal space, since both $N$-widths are equal for sorted shift-isometric spectral spaces.}

\subsection{Non half-wave symmetric functions}
Some of our results can also be extended to non half-wave symmetric functions.
However, we were only able to derive an estimate for the $N$-width and not a representation as before.
Recall from \S\ref{Subsec:Fourier} that any  $g \in L_2(\OP)$ has a unique decomposition $g=g^\even + g^\odd$ into its even and odd HWS part.

\begin{proposition} \label{prop-non-hw-symm}
    Let $N \in \mathbb{N}$, $g =g^\even + g^\odd\in L_2(\OP)$ and let the eigenvalues be ordered as $\big\{ \lambda_0^{g^\even}, \lambda_k^{g^\even}, \lambda_k^{g^\odd}: k\in\N \big\} =: \big\{ \bar\lambda^g_1\ge \bar\lambda^g_2\ge\cdots \big\}$. Then,
    \vspace*{-0.1cm}
    \begin{align*}
        \delta_N( \cU^g )^2
        \leq d_N( \cU^g )^2
        \leq  2 \sum_{\ell=N+1}^\infty \bar\lambda^g_{\ell} .
    \end{align*}
\begin{proof} 
Let $g=g^\even + g^\odd$, with $g^\even \in L_2^\even, g^\odd \in L_2^\odd$, be uniquely decomposed. Denote by $H_N^g$ the space of dimension $N$ spanned by the eigenfunctions according to the largest eigenvalues \rev{$\bar\lambda^g_1,..., \bar\lambda^g_N$. Thus, there is $M \in \N_0$ such that $ \sum_{\ell=1}^{N} \bar\lambda^g_{\ell} = \sum_{\ell=1}^{M} \lambda^{g^{\even}}_{\ell-1} + \sum_{\ell=1}^{N-M} \lambda^{g^{\odd}}_{\ell}.$} $H_N^g$ can be decomposed (not necessarily orthogonal) as $H_N^g = H_{M}^{g^\even} \rev{+} H_{N-M}^{g^\odd}$ of dimension $M\le N$ and $N-M$. Then, denoting by $P_N^g$, $P_M^{g^\even}$ and $P_{N-M}^{g^\odd}$  the orthogonal projection onto $H_N^g$, $H_{M}^{g^\even}$ and $H_{N-M}^{g^\odd}$, respectively, and recalling Remark \ref{Remark:Widths} (ii) yields
\vspace*{-0.2cm}
\begin{align*}
	& \delta_N( \cU^g )^2
		\leq d_N( \cU^g )^2
		\le \big\| g(\cdot-\mu) - P_N^g  g(\cdot-\mu) \big\|_{L_{\infty}(\cP;L_2(\Omega))}^2\\
	&\rev{ \le 2 \, \big\| (I-P_N^{g}) [g^\even(\cdot-\mu)] \, \big\|_{L_{\infty}(\cP;L_2(\Omega))}^2
		+  2 \, \big\| (I-P_N^{g})[ g^\odd(\cdot-\mu)] \, \big\|_{L_{\infty}(\cP;L_2(\Omega))}^2 } \\ 
	&  \rev{ \le 2 \, \big\| (I-P_M^{g^\even})[ g^\even(\cdot-\mu)] \, \big\|_{L_{\infty}(\cP;L_2(\Omega))}^2
		+  2 \, \big\| (I-P_{N-M}^{g^\odd})[ g^\odd(\cdot-\mu)] \, \big\|_{L_{\infty}(\cP;L_2(\Omega))}^2  }  \\
	& \rev{= 2\sum_{\ell =M+1}^\infty \lambda^{g^{\even}}_{\ell-1} + 2\sum_{\ell =N-M+1}^\infty \lambda^{g^{\odd}}_{\ell}} 
	= 2\sum_{\ell =N+1}^\infty \bar\lambda^g_{\ell},
\end{align*}
\vspace*{-0.2cm}
	which concludes the proof.
\end{proof}
\end{proposition}
\rev{Finally, by weakening} Proposition \ref{prop-non-hw-symm}, we obtain the following estimate. 
\begin{corollary}
\rev{Let $N \in \N$, then
	$\delta_{N}( \cU^g )^2
	\leq  2 \, \delta_{M}( \cU^{g^{\even}} )^2 + 2 \, \delta_{N-M}( \cU^{g^{\odd}} )^2$  
    for every $M \in \{ 0, 1, ..., N \}$.}
    \qed
\end{corollary}

\section{The effect of smoothness on the \textit{N}-width}
\label{Sec:Smoothness}
So far, we did not use any specific properties of the function $g$ modeling initial and boundary values of the linear transport equation. In this section, we shall investigate the influence of the regularity \rev{(in the sense of \S\ref{sobolev-spaces})} to the decay of the $N$-width. In particular, we use the above exact representation to derive formulae for the decay of the $N$-width.
\subsection{\rev{Sharp upper bound for} finite regularity}
\rev{For $g \in H^{r,\even}(\OP) = L_2^{\even} \cap H^{r}(\OP)$, the solution $u_{\mu}=g(\cdot-\mu) \in H^{r}(\Omega)$ is (even) periodic on $\Omega=(0,1)$, i.e., $\mathcal{U}^g \subset H^{r}(\Omega)$, recall \eqref{eq:Hr}}. For integer $r\in\N$, these spaces coincide with standard Sobolev spaces, for which Kolmogorov's classical result is known, see \cite{kolmogoroffUberBesteAnnaherung1936} and \cite[Theorem 1]{pinkusNWidthsApproximationTheory1985}.
\begin{theorem}[Kolmogorov]
	Let $r\in\N$ be an integer, then $d_N(B_{1})= (2\pi)^{-r} \,\left\lceil\tfrac{N}{2}\right\rceil^{-r}$ 
	\rev{with closed ball $B_{1}=B_{1}[H^{r}(\Omega)]$ defined in \eqref{eq:def:ball}.}\hfill\qed
\end{theorem}
This indicates that \rev{the rate} $\cO(N^{-r})$ sets a benchmark for what we can hope to achieve. \\
\rev{For the transport problem, we are going to prove the sharp upper bound for the (broken) regularity $r\in\R^+$. To show it, we first need an auxiliary lemma, which will be the main ingredient of the proof.}
\begin{lemma}\label{lem-sup-quotient}\rev{Let $r >0$, $M \in \N$, then $\quad \sup\limits_{\bc\in \overset{\searrow}{\ell}_{2,r} }
		\tfrac{
		\sum\limits_{k = M+1}^{\infty}  |c_{k}|^2
		}{
		\sum\limits_{k=1}^{\infty} k^{2r} |{c}_{k}|^2
		} = \Big( \min\limits_{m \in \N} \, \tfrac{1}{m} \sum_{k=1}^{M+m}k^{2r} \Big)^{-1} $, \\
		where $\overset{\searrow}{\ell}_{2,r}:=\{ \bc\in\ell_2(\N) \text{ decreasing}:\, \sum\limits_{k=1}^{\infty} (1+k^{2r}) |c_{k}|^2<\infty\}$.}
\begin{proof}
\rev{We start by observing that any decreasing $\bc \in \overset{\searrow}{\ell}_{2,r}$ can be represented by an sequence $\bm{a} \in \ell_2(\N)$ with 
$$ |{a}_{i}|^2 := |{c}_{i}|^2 - |{c}_{i+1}|^2 \geq 0 , \quad i \in \N,  \qquad \iff \qquad |{c}_{k}|^2 = \sum_{i=k}^{\infty} |{a}_{i}|^2 .$$
Then, we have by summation rules and fraction arithmetic
\begin{align*}
\frac{
		\sum\limits_{k = M+1}^{\infty}  |c_{k}|^2
		}{
		\sum\limits_{k=1}^{\infty} k^{2r} |{c}_{k}|^2
		} = 
		\frac{
		\sum\limits_{k = M+1}^{\infty}  \sum\limits_{i=k}^{\infty} |{a}_{i}|^2
		}{
		\sum\limits_{k=1}^{\infty} k^{2r} \sum\limits_{i=k}^{\infty} |{a}_{i}|^2
		} 
		= \frac{
		\sum\limits_{i = M+1}^{\infty} (i-M) |{a}_{i}|^2
		}{
		\sum\limits_{i = 1}^{\infty} |{a}_{i}|^2 \sum\limits_{k = 1}^{i} k^{2r}
		} \leq
		\frac{
		\sum\limits_{i = M+1}^{\infty} (i-M) |{a}_{i}|^2
		}{
		\sum\limits_{i = M+1}^{\infty} |{a}_{i}|^2 \sum\limits_{k = 1}^{i} k^{2r}
		} = 
		\frac{
		\sum\limits_{j=1}^{\infty} |{b}_{j}|^2
		}{
		\sum\limits_{j=1}^{\infty} |{b}_{j}|^2 G_M(j)
		},
\end{align*}
where the sequence $\bm{b} \in \ell_2(\N)$ in the last equality is defined by $|{b}_{j}|^2 := j \, |{a}_{M+j}|^2 \geq 0$, $j \in \N$, and $G_M(j) := \tfrac{1}{j} \sum_{k=1}^{M+j}k^{2r}$ is a discrete function. So, we simplified the problem to 
$$\quad \sup_{\bc\in \overset{\searrow}{\ell}_{2,r} }
		\frac{
		\sum\limits_{k = M+1}^{\infty}  |c_{k}|^2
		}{
		\sum\limits_{k=1}^{\infty} k^{2r} |{c}_{k}|^2
		} = \sup_{\bm{b} \in {\ell}^M_{2,r}} \frac{
		\sum\limits_{j=1}^{\infty} |{b}_{j}|^2
		}{
		\sum\limits_{j=1}^{\infty} |{b}_{j}|^2 G_M(j)
		} , $$
where ${\ell}_{2,r}^M:=\{ \bm{b} \in\ell_2(\N) :\, \sum\limits_{j=1}^{\infty} (1+G_M(j)) |b_{j}|^2<\infty\}$. We have $G_M(j) > 0$, $j \in \N$, and by 
$$ G_M(j) =  \tfrac{1}{j} \sum_{k=1}^{M+j}k^{2r} \geq \tfrac{1}{j} \sum_{k=1}^{1+j}k^{2r} > \tfrac{1}{j} \int_{0}^{j} x^{2r} dx = \frac{j^{2r}}{2r+1} \overset{(j \to \infty)}{\longrightarrow} \infty, $$
the existence of a positive minimum $\xi \in \N$, $G_M(\xi)=\min_{j \in \N} G_M(j)$. Finally, with
$$ \min_{j \in \N} G_M(j) \,  \sum\limits_{j=1}^{\infty}  |{b}_{j}|^2 =  \sum\limits_{j=1}^{\infty}  |{b}_{j}|^2 G_M(\xi)  \leq \sum\limits_{j=1}^{\infty} |{b}_{j}|^2 G_M(j) ,$$
and the minimum-generating sequence $\bm{b}^* \in {\ell}_{2,r}^M$, $|{b}^*_{j}| := \delta_{j,\xi}$, it follows that 
$$  \tfrac{1}{\min_{j \in \N} G_M(j)} = \tfrac{
		\sum\limits_{j=1}^{\infty} |{b}^*_{j}|^2
		}{
		\sum\limits_{j=1}^{\infty} |{b}^*_{j}|^2 G_M(j)} \leq  \sup_{\bm{b} \in {\ell}_{2,r}^M } \tfrac{
		\sum\limits_{j=1}^{\infty} |{b}_{j}|^2
		}{
		\sum\limits_{j=1}^{\infty} |{b}_{j}|^2 G_M(j)
		} \leq \tfrac{1}{\min_{j \in \N} G_M(j)} $$ 
		and the lemma is proven.}
\end{proof}
\end{lemma}
\begin{theorem}\label{sol-manifold-decay}
\rev{It holds for $N \in \N$ odd that
	\begin{equation} \label{Thm:sup-decay}
		\sup_{\substack{g\in H^{r,\even}, \\ \lvert g \rvert_{H^r(\OP)}\leq \sqrt{2}}} d_{N}(\mathcal{U}^g)
 = (2 \pi)^{-r} \Big( \min_{m \in \N} \, \tfrac{1}{m}  \sum_{k=1}^{ {\scriptstyle \frac{N-1}{2}}+m} k^{2r} \Big) ^{-1/2} .
	\end{equation}
	}
\begin{proof}
\rev{We first define $B_{\sqrt{2}}^{r,\HWS} := \big\{ g\in H^{r,\HWS} : \lvert g \rvert_{H^r(\OP)}\le\sqrt{2} \big\}$. Let us start by justifying that the quantity of interest is bounded for $N \geq 1$
\begin{align*}
	d_{N}(\mathcal{U}^g)^2 & \leq \dist(\Span \{ 1 \}, \cU^g)_{L_{\infty}(\cP)}^2 = \sum_{k \in \Z}^{\infty} \tfrac{|\hat{c}_{2k}(g)|^2 }{2} \leq \sum_{k \in \Z}^{\infty} (2\pi k )^{2r} \tfrac{|\hat{c}_{2k}(g)|^2 }{2} = \tfrac{\lvert g \rvert_{H^r(\OP)}^2 }{2} , \\
	 & \Longrightarrow \quad D_N := \sup_{g \in B_{\sqrt{2}}^{r,\HWS}} d_{N}(\mathcal{U}^g) \leq \sup_{g \in B_{\sqrt{2}}^{r,\HWS}} \tfrac{\lvert g \rvert_{H^r(\OP)} }{\sqrt{2}} \leq 1.
\end{align*}
On the other side, for any $ g \in B_{\sqrt{2}}^{r,\HWS}$, it also holds $g + c \in B_{\sqrt{2}}^{r,\HWS}$ for all $c \in \R$. Therefore, if we don't include constants, i.e., $\Span \{ 1 \} \not\subset H_N$, we get 
\begin{align*}
		\sup_{g \in B_{\sqrt{2}}^{r,\HWS}} \dist(H_N , \cU^g)_{L_{\infty}(\cP)}^2
		= \sup_{g \in B_{\sqrt{2}}^{r,\HWS}, \, c \in \R} \Big( c^2 + \dist(H_N , \cU^{g})_{L_{\infty}(\cP)}^2  \Big) = \infty.
	\end{align*}
	Hence, we will always include $W_0 = \Span \{ 1 \}$ and choose the remaining $N-1$ basis functions from the shift-invariant decomposition. We choose $N$ as odd and set $M=(N-1)/2 \in \N$. Then, by Theorem \ref{Theorem:Optimal-simple} and Lemma \ref{Lem:EVeven}
\begin{align*}
	d_N(\cU^g)^2 
	&= \sum_{k = N+1}^{\infty} \lambda_{k}^{g}
	= \tfrac12 \sum_{k = M+1}^{\infty}  (|\hat{c}_{\sigma(2k)}(g)|^2+|\hat{c}_{-\sigma(2k)}(g)|^2)
	= \sum_{k = M+1}^{\infty}  |\hat{c}_{\sigma(2k)}(g)|^2,
\end{align*}
where $\sigma$ again denotes an optimal sorting. Then we can express the quantity $D_N$ by
\begin{align*}
	D_N^2
	&= \sup_{g \in B_{\sqrt{2}}^{r,\HWS}} d_{N}(\mathcal{U}^g)^2 = \sup_{g\in H^{r,\even}(\OP)}
		\frac{
		2 \, d_{N}(\mathcal{U}^g)^2
		}{
		\lvert g \rvert_{H^r(\OP)}^2
		}  \\
		&= \sup_{g\in H^{r,\even}(\OP)}
		\frac{ 2\, 
		\sum\limits_{k = M+1}^{\infty}  |\hat{c}_{2\sigma(k)}(g)|^2
		}{ 2 \,
		\sum\limits_{k=1}^{\infty} (2\pi k)^{2r} |\hat{c}_{2k}(g)|^2
		} = (2\pi)^{-2r} \, \sup_{\bc\in\ell_{2,r}}
		\frac{
		\sum\limits_{k = M+1}^{\infty}  |c_{\sigma(k)}|^2
		}{
		\sum\limits_{k=1}^{\infty} {k}^{2r} |{c}_{k}|^2
		} ,
	\end{align*}
For any $\bc\in \ell_{2,r}$ and sorting of the coefficients $\sigma$ it holds by the rearrangement inequality
$$ \sum\limits_{k=1}^{\infty} k^{2r} |{c}_{k}|^2 \geq \sum\limits_{k=1}^{\infty} k^{2r} |{c}_{\sigma(k)}|^2 \quad \iff \quad \frac{
		\sum\limits_{k = M+1}^{\infty}  |c_{\sigma(k)}|^2
		}{
		\sum\limits_{k=1}^{\infty} k^{2r} |{c}_{k}|^2
		} \leq \frac{
		\sum\limits_{k = M+1}^{\infty}  |c_{\sigma(k)}|^2
		}{
		\sum\limits_{k=1}^{\infty} k^{2r} |{c}_{\sigma(k)}|^2
		},  $$
which means we can restrict ourselves to $\sigma\equiv\text{Id}$, i.e., already sorted sequences. \\
Putting all together gives 
$$  D_N^2 = (2\pi)^{-2r} \sup_{\bc\in\ell_{2,r}}
		\tfrac{
		\sum\limits_{k = M+1}^{\infty}  |c_{k}|^2
		}{
		\sum\limits_{k=0}^{\infty} k^{2r} |{c}_{\sigma(k)}|^2
		} = (2\pi)^{-2r} \sup_{\bc\in \overset{\searrow}{\ell}_{2,r} }
		\tfrac{
		\sum\limits_{k = M+1}^{\infty}  |c_{k}|^2
		}{
		\sum\limits_{k=0}^{\infty} k^{2r} |{c}_{k}|^2,
		} $$  
where the statement of the theorem follows by Lemma \ref{lem-sup-quotient}.}
\end{proof}
\end{theorem}
\begin{remark}\label{rem-sup-decay}
\rev{\begin{compactenum}[(a)]
	\item For general $g \in H^{r,\even}$, with Theorem \ref{sol-manifold-decay} we have
	$$d_{N}(\cU^g) \leq \lvert g \rvert_{H^r(\Omega)} \, (2 \pi)^{-r} \Big( \min_{m \in \N} \, \tfrac{1}{m} \sum_{k=1}^{ {\scriptstyle \frac{N-1}{2}}+m} k^{2r} \Big) ^{-1/2}  ,$$ 
	which might be substantially smaller as we picked really a special case as the supremum. 
   \item Trivially, on $H^{r,\even}$ without further assumptions, there can't be a lower bound since even $d_{N}(\cU^g)^2 \equiv 0$ for $g \equiv 0 \in H^{r,\eo}$.
   	\item If we use non-optimal spaces of Fourier modes just ordered by frequency, i.e., $X_N = \bigoplus_{k=0}^{(N-1)/2} W^{\even}_k$, then we get for $N$ odd
   	$$  \sup_{g \in B^{r,\even}_{\sqrt{2}}} \dist(X_N , \cU^g)_{L_{\infty}(\cP)} =  (2\pi)^{-r} \, \big( \tfrac{N+1}{2} \big)^{-r} .$$
   	\begin{proof}
   	We follow the steps of the proof in Theorem \ref{sol-manifold-decay}, but the expression becomes much simpler if we don't have a sorting. We directly have\\ $\sup\limits_{\bc\in \ell_{2,r} }
		\tfrac{
		\sum\limits_{k = M+1}^{\infty}  |c_{k}|^2
		}{
		\sum\limits_{k=1}^{\infty} k^{2r} |{c}_{k}|^2
		} = (M+1)^{-2r}$ instead of Lemma \ref{lem-sup-quotient}. 
   	\end{proof}
   	\item For the odd HWS counterpart, we have \\
    $ \sup_{g \in B^{r,\odd}_{\sqrt{2}}} \delta_{N}(\cU^g) = (\pi)^{-r} \Big( \min_{m \in \N} \, \tfrac{1}{m}  \sum_{k=1}^{ {\scriptstyle \frac{N}{2}}+m} (2k-1)^{2r} \Big) ^{-1/2}$.
   	\begin{proof}
   	Follows from exactly the same arguments as in the proof of Theorem \ref{sol-manifold-decay}, except with an odd-frequency basis with no constant function.
   	\end{proof}
   	\item By losing a factor of 2, the same can be applied to non-HWS functions using Proposition \ref{prop-non-hw-symm}. \label{rem-sup-decay-f}
\end{compactenum}
}
\end{remark}
$\text{ }$\\
However, if we know a more specific form of the Fourier coefficients, we can get an lower estimate. Further, we also get a reverse implication.

\subsection{\rev{Exact rate by eigenvalue decay}}\label{Sec:Exact}
\smallskip
\subsubsection*{Exact decay of some piecewise functions}
So far, we derived \emph{upper bounds} for the decay of the $N$-widths for $g \in H^r(\OP)$.
In this section, we investigate exact formulae for the decay as well as \emph{lower bounds}.
First, we are going to investigate specific functions of known regularity, for which lower bounds can be proven.
These functions are such that $g^{(m)}=\operatorname{sgn}(\cdot)$, $m\in\N$, i.e., the jump function considered in Example~\ref{exam-jump-L2}.
This is done in a recursive manner by setting $g_0:=\operatorname{sgn}(\cdot)$.
Then, for $m = 1,2,3,...$, we set
\begin{equation}\label{define_fn}
	g_{m}(y):=\int_{0}^{y} g_{m-1}(x)\, dx -\tfrac{1}{2} \int_{0}^{1} g_{m-1}(x)\, dx, \quad y \in [-1,1],
\end{equation}
which results in piecewise polynomials,
\begin{align*}
    g_{1}(y) = \begin{cases}  - y - \tfrac{1}{2} , \quad y \in [-1,0), \\ 
                   + y - \tfrac{1}{2} , \quad y \in [0,+1] ,
                   \end{cases}   
    \quad  g_{2}(y) 
            = \begin{cases}  - \tfrac{1}{2} y^2 - \tfrac{1}{2}y , \quad y \in [-1,0),  \\ 
                +  \tfrac{1}{2} y^2 - \tfrac{1}{2}y , \quad y \in [0,+1].\end{cases}
\end{align*}

\begin{lemma} \label{fn_integral_of_jump}
    Let $ m \in\N_0$ and $g_m$ as defined in \eqref{define_fn}.
    Then, $ g_m \in H^{ r(m) -\varepsilon, \odd}$ for $r(m)=m+1/2$ and all $\varepsilon>0$ but $g_m \not\in H^{r(m), \odd}$, $g_m^{(m)} = g_{0} =  \operatorname{sgn}(\cdot)$ and
    \vspace*{-0.7em}
    \begin{align}\label{eq:fn_integral_of_jump}
        \| Q^\odd_{k} g_{ m} \|_{L_2(\Omega)} = \sqrt{8} ((2k-1)\pi)^{-({ m }+1)} = \rev{\sqrt{\lambda^g_{2k} + \lambda^g_{2k-1}}}.
    \end{align}
\begin{proof}
	It is well-known (e.g.\ by considering the Fourier expansion), that $\operatorname{sgn}(\cdot)\in H^{1/2-\varepsilon}(\OP)$ for all $\varepsilon>0$, but  $\operatorname{sgn}(\cdot)\not\in H^{1/2}(\OP)$. Then, the recursive definition immediately implies the statement about the regularity. Moreover, $g_{0}\in L_2^\odd$. Assuming that $g_{{m}-1}\in L_2^\odd$ for ${m}\in\N$, we get for any 
    $-1<x<0$
    \begin{align*}
    g_{m}(x) + g_{m}(x+1) 
    &= \int_{0}^{x} g_{{m}-1}(y)\, dy + \int_{0}^{x+1} g_{{m}-1}(y)\, dy -2 \, \tfrac{1}{2} \int_{0}^{1} g_{{m}-1}(y)\, dy\\
    &= \int_{0}^{x} g_{{m}-1}(y)\, dy + \int_{1}^{x+1} g_{{m}-1}(y)\, dy = 0,
    \end{align*}
    by Lemma \ref{Lem:HWS}, so that $g_{m}\in L_2^\odd$ by induction.
    Finally, Example~\ref{exam-jump-L2} yields
    \begin{align*}
        \| Q_{k}^\odd g_{m} \|_{L_2(\Omega)}
        \kern-3pt=\kern-2pt ((2k-1)\pi)^{-{m}} \| Q_{k}^\odd g_{m}^{({m})} \|_{L_2(\Omega)}
        \kern-3pt=\kern-2pt \sqrt{8} ((2k-1)\pi)^{-({m}+1)},
    \end{align*}
    which proves \eqref{eq:fn_integral_of_jump} and finishes the proof.
\end{proof}
\end{lemma}

\begin{theorem} \label{Thm-decay-lower-bound}
    Let $ m \in\N_0$ and $g_m$ as defined in \eqref{define_fn}. Then,
   \begin{align*}
        \delta_{N}(\cU^g)^2
        &= \tfrac{4}{\pi^{2m+2}} \! \sum_{\ell= \lceil (N+1)/2 \rceil }^{\infty}  (2 \ell -1)^{-2{m}-2} ,
    \end{align*}
    and, for even $N$, we have $\delta_{N}(\cU^g) = d_{N}(\cU^g)$. Moreover
       $$
 \delta_{N}(\mathcal{U}^g) \; \in \; \tfrac{ \sqrt{2} }{\sqrt{2m+1} \, \pi^{m+1}} \, \left[ \big( N+1 \big)^{-m-\frac{1}{2}}, \, \big( N-2 \big)^{-m-\frac{1}{2}} \right]
$$
for all $N\in\N$. Thus $\delta_{N}(\cU^g)\cong N^{-r(m)}$ in terms of $r(m)$.
\begin{proof}
First, by \eqref{eq:fn_integral_of_jump} the trivial sorting $\sigma \equiv \operatorname{Id}$ is the optimal sorting from \eqref{opt-sorting} and therefore we get by Theorem \ref{Theorem:Optimal-simple} the representation 
$$ \delta_{N}(\cU^g)^2
        = 4\sum_{k=N+1}^{\infty}  \Big( (2 \lfloor \tfrac{k+1}{2} \rfloor -1)  \pi \Big)^{-2{m}-2} = \tfrac{4}{\pi^{2m+2}} \! \sum_{\ell= \lceil (N+1)/2 \rceil }^{\infty}  (2 \ell -1)^{-2{m}-2} , $$ 
        which is equal to $d_N(\cU^g)$ for even $N$. To prove the bounds, we deduce that on the one hand, 
	\begin{align*}
  \tfrac{\delta_{N}(\mathcal{U}^g)^2 \, \pi^{2m+2}}{4} = \hspace{-12pt} \sum_{\ell=\lceil (N+1)/2 \rceil}^{\infty} \hspace{-6pt} \left( \tfrac{1}{2\ell-1} \right)^{2m+2} \leq \int_{\lceil (N+1)/2 \rceil -1}^{\infty} \left( \tfrac{1}{2x-1} \right)^{2m+2} \, dx = \tfrac{\big(2 \lceil (N + 1)/2 \rceil -3 \big)^{-(2m+1)}}{2(2m+1)},
 \end{align*}
 and on the other hand
    \begin{align*}
  \tfrac{\delta_{N}(\mathcal{U}^g)^2 \, \pi^{2m+2}}{4} = \hspace{-12pt} \sum_{\ell=\lceil (N+1)/2 \rceil}^{\infty} \hspace{-6pt} \left( \tfrac{1}{2\ell-1} \right)^{2m+2} \geq \int_{\lceil (N+1)/2 \rceil }^{\infty} \left( \tfrac{1}{2x-1} \right)^{2m+2} \, dx = \tfrac{\big(2 \lceil (N + 1)/2 \rceil - 1 \big)^{-(2m+1)}}{2(2m+1)}.
 \end{align*}
 Then by
 $$ \left[ \big(2 \lceil (N + 1)/2 \rceil - 1 \big)^{-1}, \, \big(2 \lceil (N + 1)/2 \rceil - 3 \big)^{-1} \right] \subset \left[ \big( N+1 \big)^{-1}, \, \big( N-2 \big)^{-1} \right] $$
we get the desired bounds. 
\end{proof}
\end{theorem}

\rev{Let us capture this relationship between eigenvalue decay and $N$-width decay more generally.}

\subsubsection*{\rev{Exact rate by eigenvalue decay}}
\rev{In section \ref{section3}, Lemma \ref{Eigs-are-optimal}, we proved for any problem we can express
\begin{align*}
\delta_N(\cU^g)^2 =  \sum_{k=N+1}^{\infty} \lambda_{k}^g \; ,
\end{align*}
which also equals $d_N(\cU^g)^2$ for the case of shift-isometric spaces. We aim for an expression of the type $\delta_N(\cU^g)^2 = A(N),$ where $A$ is some algebraic term. But generally, the sequence $(\lambda_{k}^g)_{k \in \N} \in \ell_1$ can have any form at first. However, if we can express either the eigenvalues or the $N$-width algebraically, we can express both algebraically.
}

\begin{lemma} \label{lem:discrete-ftc}
\rev{
Let $f:[1,\infty) \to \R_{+}$ be an decreasing and improper integrable function with antiderivative $F:[1,\infty) \to \R_{-}$ defined by $F(y):= - \int_{y}^{\infty} f(x) dx, \; y \geq 1$.\\
Then we arrive at a discrete version of the fundamental theorem of calculus.
\begin{align*}
&(a)& \quad \lambda_{n}^g &= f(n), \; n \in \N  &\quad \Longrightarrow \quad& -F(n+1) \leq \delta_{n}(\cU^g)^2 \leq -F(n), \; n \in \N  \\
&(b)& \quad \delta_{n}(\cU^g)^2 &= -F(n), \; n \in \N  &\quad \Longrightarrow \quad& f(n) \leq \lambda_{n}^g \leq f(n-1), \; n \in \N_{\geq 2}
\end{align*}
We can relax the first statement to better fit our objectives. Let $C \geq c \geq 0$, then
\begin{align*}
&(a') \quad c f(n) \leq \lambda_{n}^g \leq C f(n), \; n \in \mathbb{N} \quad \Longrightarrow \; -c \, F(n+1) \leq \delta_{n}(\mathcal{U}^g)^2 \leq -C \, F(n), \; n \in \mathbb{N}.
\end{align*}
}
\begin{proof}
\rev{
"$(a')$" For $C \geq 0$ and $n \in \N$ \footnote{Using the Euler–Maclaurin formula, the estimate can be done even more tight. But we didn't want to delve deeper into these technicalities; this should suffice for the asymptotic behavior.}
\begin{align*}
\delta_n(\cU^g)^2 \leq \sum_{k=n+1}^{\infty} C \, f(k) \leq C \sum_{k=n+1}^{\infty}  \int_{k-1}^{k} f(x) dx = C \int_{n}^{\infty} f(x) dx = - C \, F(n).
\end{align*}
Similarly, for $c \geq 0$ and $n \in \N$
\begin{align*}
\delta_n(\cU^g)^2 \geq \sum_{k=n+1}^{\infty} c \, f(k) \geq c \sum_{k=n+1}^{\infty}  \int_{k}^{k+1} f(x) dx = c \int_{n+1}^{\infty} f(x) dx = - c \, F(n+1).
\end{align*}
"$(b)$" We have for $n \in \N_{\geq 2}$ 
\begin{align*}
\lambda_{n}^g = \delta_{n-1}(\cU^g)^2 - \delta_n(\cU^g)^2 = F(n) - F(n-1) = \int_{n-1}^n f(x) dx \quad 
\begin{array}{ll}
\leq & f(n-1), \\
\geq & f(n).
\end{array}
\end{align*}
}
\end{proof}
\end{lemma}
\rev{
\begin{remark} \label{rem-Hrminuseps}
\begin{compactenum}[(a)]
\item Defining $g \in L_2^{\even}$ by $ \tfrac12 \lvert \hat{c}_{2k}(g) \rvert^2:= k^{-(2r+1)} = \lambda^g_{2k-1} = \lambda^g_{2k}$ results in a function $g\in H^{r-\varepsilon,\even}$ for all $\varepsilon>0$, but $g\not\in H^{r,\even}$ with $\delta_N(\cU^g)$ being strictly smaller than $c \, N^{-r}$.
	\item Vice versa, using Lemma \ref{lem:discrete-ftc}, our analysis shows that $d_N(\cU^g) = c\, N^{-r}$ implies the corresponding Sobolev regularity $r-\varepsilon$ for all $\varepsilon>0$ of $g$.
\item Defining $g \in L_2^{\even}$ by $ \tfrac12 \lvert \hat{c}_{2k}(g) \rvert^2:= k^{-(2r+1)} \log(k)^{-2} = \lambda^g_{2k-1} = \lambda^g_{2k}$ results in a function $g\in H^{r,\even}$, but  $g\not\in H^{r+\delta,\even}$ for all $\delta >0$, with $\delta_N(\cU^g) \approx C N^{-r} \log(N)^{-1}$ being strictly smaller than $c \, N^{-r}$.
\end{compactenum}
\end{remark}
}

\subsection{\rev{Limit case}: Infinite regularity}
\rev{In this subsection we deal with the limit case $g \in C^{\infty}(\OP)$. We address the question of whether $g$'s infinite differentiability implies exponential decay. \\
One can trivially always express $\delta_N(\mathcal{U}^g) = e^{-\gamma(N)}$ with an increasing sequence $\gamma(n) := -\log(\delta_n(\mathcal{U}^g) ),$ $n \in \N$.\\ 
Then directly from $g \in C^{\infty }(\OP) \subset H^{r}(\OP)$, we get by Remark \ref{rem-sup-decay} that $\delta_N(\mathcal{U}^g) < C_r \, N^{-r} = C_r \, e^{-r \, \log(N)}$ for all $r>0$. This means, $\gamma(n)$ grows faster than linear in the logarithm, i.e., for fixed $a,b >0$ it is 
$$\gamma(n) > a + b \, \log(n)$$ 
for $n$ sufficiently large.
\begin{example}
One example is $g \in C^{\infty}(\OP)$ with $\lambda_k^g := \tfrac{2 \log(k)}{k} k^{- \log(k)}$, $k \in \N$, which results with Lemma \ref{lem:discrete-ftc} in a decay of $$e^{- \log(N+1)^2} < \delta_N(\mathcal{U}^g) < e^{- \log(N)^2} . $$
\end{example}
But can this already be understood as exponential? The decay is faster than polynomial $C_r \, N^{-r}$, but slower than exponential $ e^{-N}$. To really get an exponential decay in the classical sense, i.e., $\gamma(n) = \operatorname{linear}$, we already need exponential decay of the coefficients, since we get straight from Lemma \ref{lem:discrete-ftc} the following equivalence.
\begin{corollary}
We get for $C > 0$ and a base $ \omega > 1$ that
\begin{align*}
\lambda_k^g \leq C \, \omega^{-2k } , \, \; k \in \N \hspace{0.5em} \Rightarrow \; \delta_N(\mathcal{U}^g) \leq \sqrt{ \tfrac{C}{2 \, \log(\omega) }} \, \omega^{- N} , \, \; N \in \N \hspace{0.5em} \Rightarrow \; \lambda_k^g \leq ( C \, \omega^2) \, \omega^{-2k } , \, \; k \in \N_{\geq 2} .
\end{align*}
\end{corollary}
$\text{ }$ \\
Next we give a proposition to show that a holomorphic function indeed implies exponential decay.\\}

\begin{proposition} 
    Let $g \in C^{\infty}(\OP) $ with a complex analytic extension $\bar{g}:\cB_1(0)\to \C$\footnote{We denote the ball of radius $\varrho>0$ around $z\in\C$ by $\cB_{\varrho}(z):=\{ y\in\C: |y-z| \leq \varrho \}$.}, $\bar{g}_{|\OP} = g$, and a constant $C^g :=\max _{z \in \mathcal{B}_{1}(0)}\lvert \, \bar{g}(z) \, \rvert > 0$. Then,
    $\delta_{N}( \cU^g )
        \leq d_{N}( \cU^g )
        <  (C^g K) \, \tau^{-N}$, where $\tau =\pi \, e^{7/8} \approx 7.536$ and $K= e^{9/8} 2^{-1/2} \approx 2.178$.
\begin{proof}
  We consider the holomorphic extension $\bar{g}$ in $ \cB_1(0) \supset (-1,1)$. Cauchy's integral formula gives $| \bar{g}^{(r)} (x) | \leq r!\, 1^r \max_{z \in \partial \cB_{1}(x)} |\bar{g}(z)| \leq C^g \, r!$, \cite{MR924157}. By half-wave symmetry and Sobolev norms, we get $\lvert g \rvert_{H^r(\Omega)} = \|g^{(r)}\|_{L_2(\Omega)} = 2^{-1/2} \, \|g^{(r)}\|_{L_2(\OP)} \leq C^g \,r!\, 2^{-1/2}$. Then, 
  \begin{equation*}
  d_{N}( \cU^g )
       < \lvert \, g \, \rvert_{H^r(\Omega)} \pi^{-r} \, N^{-r} \leq C^g \, r!\, 2^{-1/2} \pi^{-r} \, N^{-r} 
  \end{equation*}
   follows by Remark~\ref{rem-sup-decay}~\ref{rem-sup-decay-f}. Using Stirling's approximation $N! \leq e \sqrt{N} (N/e)^N $ for the factorial and setting $r=N$ gives $d_{N}( \cU^g ) < C^g \, e \, \sqrt{N/2} \, (\pi e)^{-N}$. Finally, $\sqrt{N} < (e^{1/8})^N e^{1/8}$ gives the desired result.
\end{proof}
\end{proposition}%
\rev{We end this subsection with a small final point, namely that if the norm of the $r$-th (weak) derivative itself can be bounded by just a power function, we even get a \textit{finite decay}.
\begin{remark}
Let $g \in C^{\infty}(\OP)$. If there exists $D >0$ sucht that $ \frac{\| g^{(r)} \|_{L_2(\Omega)}}{\| g \|_{L_2(\Omega)}} = \frac{\lvert g \rvert_{H^r(\Omega)}}{\| g \|_{L_2(\Omega)}}  \leq D^r ,$ for all $r>0$, then it holds $ \quad d_N(\mathcal{U}^g) = 0, \quad \text{for } N >  D/\pi$. 
\begin{proof} $d_N(\mathcal{U}^g) < 2 \lvert g \rvert_{H^r(\Omega)} \, \pi^{-r} \, N^{-r} \leq 2 \| g \|_{L_2(\Omega)} \, \Big( \tfrac{ D }{ \pi N } \Big)^r  \; \xrightarrow[(r \to \infty)]{} \; 0, \quad \text{ for } N > D/\pi.$
\end{proof}
\end{remark}
}

\section{Numerical experiments} \label{Sec:experiments}
We are now going to report results of some of our numerical experiments highlighting different quantitative aspects of our previous theoretical investigations. More details and the code can be found in \cite{florianarbesFlabowskiN_widths_for_transportNumerical2024}.  

\subsection{Numerical approximation of the \textit{N}-width} \label{numerical_error}
It is clear that we cannot compute $\delta_N(\cU^g)$ or $d_N(\cU^g)$ exactly, at least in general.
Even for a given linear approximation space $V_N$, the distance of $\cU^g$ to $V_N$ amounts to computing an integral over $\cP$ in the $L_2$-case or the determination of a supremum in the $L_\infty$-framework. Both would only be possible exactly, if we had a formula for the error $u_\mu-P_N u_\mu$ at hand.

Otherwise, we need a discretization in space $\Omega$ and for the parameter set $\cP$ in such a way that the resulting numerical approximation is sufficiently accurate.
In space, we fix a number $n_x\in\N$ of uniformly spaced quadrature or sampling points $x_i$, $i=1,...,n_x$ (we choose $n_x=2500$), by setting $\Delta x := 1/n_x$ and $x_i:= (2i-1)/2\, \Delta x$.
We collect these points in a vector $\bx:=(x_1,...,x_{n_x})^\top\in\R^{n_x}$.
We proceed in a similar manner for $\cP$ by choosing $n_\mu\in\N$ of uniformly spaced points $\mu_j$, $j=1,...,n_\mu$ (we choose $n_\mu=2500$) by $\Delta\mu := 1/n_\mu$ and $\mu_j:= (2j-1)/2\, \Delta\mu$, $\bmu:=(\mu_1,...,\mu_{n_\mu})^\top\in\R^{n_\mu}$.
This corresponds to the midpoint rule for numerical integration.

\subsubsection*{Proper Orthogonal Decomposition (POD)}
For some given parameter value $\mu_j$ and a given function $g$, we determine $X_{i,j}:= g(x_i - \mu_j) = u_{\mu_j}(x_i)$ as a \enquote{snapshot} of \eqref{eq:functiontobeapproximated}.
These values are collected in the \emph{snapshot matrix} $\bX:=(X_{i,j})_{i=1,...,n_x; j=1,...,n_\mu}\in\R^{n_x\times n_\mu}$.

For the $L_2$-width, we perform a singular value decomposition (SVD) $\bX  =  \bU \bSigma\bV^\top$, which is then truncated to dimension $N\in\N$ in order to obtain a reduced basis, which corresponds to the \emph{Proper Orthogonal Decomposition (POD)}. Since it is known that the POD is the best approximation w.r.t.\ $L_2(\cP)$, we get the optimal spaces $V_N$ in the $L_2(\cP)$-sense. The error and thus the $\delta_N$-width can be computed from the singular values: $\delta_N( \cU^g )^2 \approx \sum_{k=N+1}^{\min(n_\mu, n_x)} \sigma_k^2$.

\subsubsection*{Optimal spaces}
In some cases, we have constructed (optimal) spaces $V_N$, i.e., we know an ON-basis for $V_N$.
In such a case, the SVD is given through the projection terms, i.e. no eigenvalue decomposition needs to be performed.
In case we know that $d_N( \cU^g ) = \delta_N( \cU^g )$, no further computations are needed.
In case they are not equal, we can proceed with the basis $V_N$ in order to compute the approximation error as described now in detail.

\subsubsection*{Computation of the distance/error}
In order to determine the distance of $V_N$ to the set $\cU^g$ of solutions, let $V_N=\Span\{ \psi_1,...,\psi_N\}$ for some ON-basis functions $\psi_\ell$, $\ell=1,...,N$. 
The best approximation of some function $u_\mu$ onto $V_N$ is the orthogonal projection, i.e., $P_N(u_\mu) = \sum_{\ell=1}^N \langle u_\mu, \psi_\ell\rangle_{L_2(\Omega)}\psi_\ell$.
The inner products are approximated by the midpoint rule, i.e.,
$\langle u_\mu, \psi_\ell\rangle_{L_2(\Omega)}
    \approx \tfrac1{n_x} \sum_{i'=1}^{n_x} u_\mu(x_{i'})\,  \psi_\ell(x_{i'})$,
so that by orthonormality
\begin{align*}
    \inf_{\tilde v_N\in V_N} \| u_{\mu_j} - \tilde v_N\|_{L_2(\Omega)}^2
     & = \| u_{\mu_j} - P_N(u_{\mu_j})\|_{L_2(\Omega)}^2                                                              \\
     & \kern-100pt\approx \tfrac1{n_x} \sum_{i=1}^{n_x} \big( u_{\mu_j}(x_{i}) \kern-2pt-\kern-2pt (P_N(u_{\mu_j}))(x_{i})\big)^2
      = ... \approx \tfrac1{n_x} \sum_{i=1}^{n_x} \left(X_{i,j} \kern-2pt-\kern-2pt \tfrac1{n_x} (\bPsi_N \bPsi_N^\top\, \bX)_{i,j} \right)^2,
\end{align*}
where $\bPsi_N=(\psi_\ell(x_i))_{i,\ell}\in\R^{n_x\times N}$.
Then, $\dist(V_N, \cU^g)_{L_\infty(\cP; L_2(\Omega))}$ is approximated by taking the maximum over $j=1,...,n_\mu$ of the latter quantity and then the square root.
As for the $L_2$-distance $\dist(V_N, \cU^g)_{L_2(\cP; L_2(\Omega))}$, by \eqref{deltaN_2}
\begin{align*}
    \dist(V_N,\cU^g)_{L_2(\cP; L_2(\Omega))}^2
     & = \| u_\mu - P_{N}u_\mu \|_{L_{2}(\cP;L_2(\Omega))}^2
    = \int_\cP \| u_\mu - P_{N}u_\mu \|_{L_2(\Omega)}^2\, d\mu                                          \\
     & \kern-100pt\approx \tfrac1{n_{\mu}} \sum_{j=1}^{n_\mu} \| u_{\mu_j} - P_{N}u_{\mu_j} \|_{L_2(\Omega)}^2
    \approx \tfrac1{n_{\mu}} \tfrac1{n_{x}} \sum_{j=1}^{n_\mu}
    \sum_{i=1}^{n_x} \left(X_{i,j} - \tfrac1{n_x} (\bPsi_N \bPsi_N^\top \bX)_{i,j} \right)^2 .
\end{align*}

\subsection{Error bound for a jump discontinuity} 
\label{heaviside}
\begin{figure}[!htb]
    \centering
    \begin{subfigure}{.5\textwidth}
        \centering
        \includegraphics[width=1\linewidth]{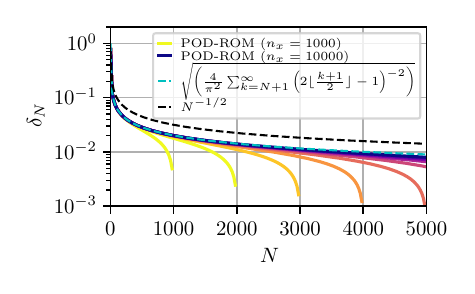}
        \caption{ POD vs.\ analytic error decay.}
        \label{fig:error_heaviside}
    \end{subfigure}%
    \begin{subfigure}{.5\textwidth}
        \centering
        \includegraphics[width=1\linewidth]{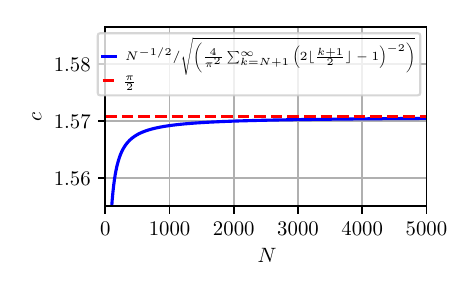}
        \caption{Ratio of exact form to asymptotic rate.}
        \label{fig:estimate_HS}
    \end{subfigure}
    \caption{Kolmorogov $N$-width $d_N( \cU^g )$ for a discontinuous function -- comparison of POD, exact form of $d_N( \cU^g )$ and known asymptotic rate.}
    \label{fig:error_bounds}
\end{figure}

For a discontinuous function, it is known that $d_N( \cU^g )\le c N^{-1/2}$ \cite{ohlbergerReducedBasisMethods2016c}, the novel exact representation is given in  \eqref{equ:exact_decay_HS}. We compare a reduced order model determined by POD  with the exact rate, which allows us to numerically  investigate the difference between the asymptotic order $N^{-1/2}$ and the exact rate.
The results are shown in Figure~\ref{fig:error_bounds}.
In the graph on the left, we show the decay for different sizes of $n_x$, i.e., various numbers of the original snapshots to build the POD (shown in different colors).
Instead of computing the SVD, we use the basis functions $\phi^\odd_{k}$ and $\psi^\odd_{k}$ defined in Lemma~\ref{La:Wodd}, as they are known to be optimal.
Numerical results confirmed that the POD basis vectors are in fact identical with the analytical basis vectors up to a seemingly random phase shift and a tolerance for numerical precision.

As we see, they asymptotically reach the exact representation shown in cyan.
This also confirms the known fact that POD is optimal w.r.t.\ the $L_2$-width.
We also show the asymptotic order $N^{-1/2}$ in black.

The formula for the exact rate cannot immediately be re-interpreted as a simple asymptotic w.r.t.\ $N$.
To this end, on the right-hand side of Figure~\ref{fig:estimate_HS} we plot the ratio of $N^{-1/2}$ and the exact form and see that it reaches $\tfrac\pi{2}$, which is interesting at least for two reasons: (i) the asymptotic rate $N^{-1/2}$ is sharp with a multiple factor of $\tfrac\pi{2}$; (ii) the exact formula has an asymptotic behavior as $N^{-1/2}$.
%

\subsection{Smooth steep functions}
\label{SubSec:steep}
\begin{figure}[!htb]
    \centering
    \captionsetup[subfigure]{width=0.9\linewidth}
    \begin{subfigure}{.5\textwidth}\centering
        \includegraphics[width=1\linewidth]{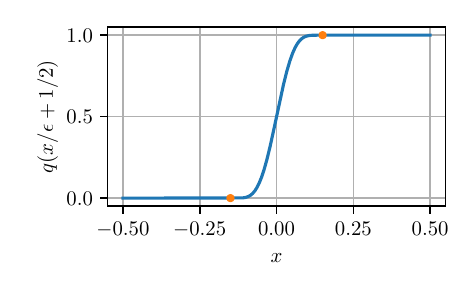}
        \caption{Smooth function $q_5: (-\tfrac12,\tfrac12)\to\R$,\\$x=\pm \varepsilon/2$ marked in orange.} 
        \label{fig:q}
    \end{subfigure}%
    \begin{subfigure}{.5\textwidth}\centering
        \includegraphics[width=1\linewidth]{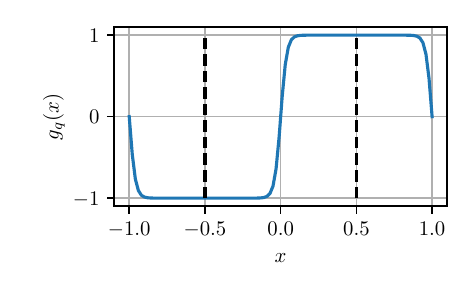}
        \caption{Resulting odd HWS $g=g_q$ with junctions at dashed lines.}
        \label{fig:g}
    \end{subfigure}
    \caption{Construction of an odd HWS initial condition from a smooth ramp.}
        \label{fig:half_wave_odd_symmetric_initial_condition}
\end{figure}

We are now considering smooth functions which are \enquote{close} to a jump in the sense that they have one or more steep ramps.
To this end, we construct an odd half-wave symmetric function, shown in Figure~\ref{fig:g}.
The starting point is some smooth odd-symmetric function $q$ on the interval $(-\tfrac12,\tfrac12)$ (see Figure~\ref{fig:q}).
Based upon this, we define the odd HWS function $g=g_q$ by
\begin{align}
    g(x) :=
    \begin{cases}
        1-2q(x+1), & -1 < x  \le -\tfrac12,      \\
        2q(x)-1,   & -\tfrac12 < x \le \tfrac12, \\
        1-2q(x-1), & \tfrac12 < x < 1.           \\
    \end{cases}
    \label{equ:half_wave_odd_symmetric_initial_condition}
\end{align}

Following this idea, we can derive functions with arbitrary smoothness and arbitrarily steep ramps in order to be able to numerically investigate the dependence of the decay rate of $d_N( \cU^g )$ on the regularity and the shape of the function.
To this end, we construct a whole family $\{ q_m\}_{m\in\N_0}$ such that $q_m\in C^m$ but $q_m\not\in C^{s
m+1}$ (so that $m$ is the exact degree of regularity of $q_m$).
We show an example of such functions $q_0,...,q_5$ in \eqref{equ:polynomials}.
Starting from a linear function $q_0$, we successively increase the polynomial degree.
A parameter $\varepsilon$ is used to control the steepness of the ramp.\footnote{We give all details for the sake of reproducible research.} Then, we get
\vspace*{-0.2cm}
\begin{subequations}\label{equ:polynomials}
    \begin{align}
        q_0(x/\varepsilon+1/2) & :=  x, \label{equ:polynomials:0}                       \\
        q_1(x/\varepsilon+1/2) & := -2 x^3 + 3 x^2,                                     \\
        q_2(x/\varepsilon+1/2) & := 6 x^5 - 15 x^4 + 10 x^3,                            \\
        q_3(x/\varepsilon+1/2) & := -20 x^7 + 70 x^6 - 84 x^5 + 35 x^4,                 \\
        q_4(x/\varepsilon+1/2) & := 70 x^9 - 315 x^8 + 540 x^7 - 420 x^6 + 126 x^5,     \\
        q_5(x/\varepsilon+1/2) & := -252 x^{11} + 1386 x^{10} - 3080 x^9 + 3465 x^8 - 1980 x^7 + 462 x^6,
    \end{align}
\end{subequations}
with the ramp being between $x=-\varepsilon/2$ and $x=\varepsilon/2$, the junctions are marked in Figure~\ref{fig:q}. Outside the ramp, $q_m \equiv 0$ and $q_m \equiv 1$ respectively. As an example of a $C^{\infty}$-function, we use the sigmoid function $q_{\infty,5}$ defined recursively as $q_{\infty,k+1}(x):=\sin \left ( \tfrac{\pi}{2} q_{\infty,k}(x) \right )$ with $q_{\infty,0}(x):=\frac{2(x-\mu)}{\varepsilon m}$, $m=\tfrac{\pi^k}{2^k}$ with the smooth limit $q_{\infty}$ having the property $q_{\infty}\equiv-1$ for $x < \mu-m\tfrac{\varepsilon}{2}$ as well as $q_{\infty}\equiv1$ for $x > \mu+m \tfrac{\varepsilon}{2}$.

\begin{figure}[!htb]\centering
    \captionsetup[subfigure]{width=0.9\linewidth}
    \begin{subfigure}{.5\textwidth}\centering
        \includegraphics[width=1\linewidth]{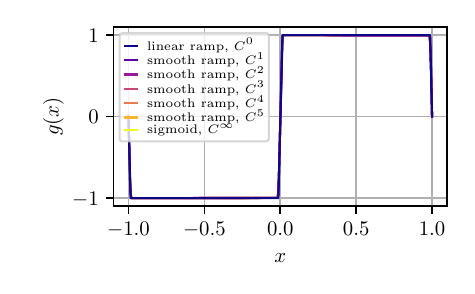}
        \caption{Ramp functions on $\OP$.\label{fig:ramps}}
    \end{subfigure}%
    \begin{subfigure}{.5\textwidth}\centering
        \includegraphics[width=1\linewidth]{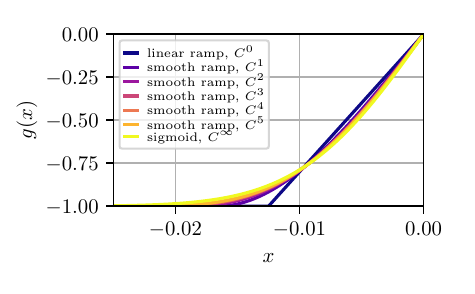}
        \caption{Zoom into subinterval $[-0.025,0]$.\label{fig:ramps2}}
    \end{subfigure}
    \caption{Ramp functions with varying smoothness $C^m$.}
    \label{fig:ramps_all}
\end{figure}

We can continue this process to obtain a $C^\infty$-function, but do not go into details.
In order to get a meaningful comparison for the dependency of the $N$-width in terms of smoothness, we will use $\varepsilon$ for such a fine-tuning.
The aim is that all functions $q_m$ should feature a similar steep jump from $0$ to $1$, but differ in their regularity, which of course causes different shapes of the functions, see Figure~\ref{fig:ramps_all}.
Hence, we fit each resulting $g_{q_m}$ to $g_{q_0}$ and choose $\varepsilon$ as the parameter resulting in the best fit.
We indicate the resulting values for $\varepsilon$ in Table~\ref{table:eps}.
\begin{table}[!htb]
    \begin{tabular}{ | m{1.4cm} | m{1.2cm} | m{1.2cm} | m{1.2cm} | m{1.2cm} | m{1.2cm} | m{1.2cm} | }
        \hline
        regularity & $C^0$ & $C^1$   & $C^2$   & $C^3$   & $C^4$   & $C^5$   \\
        \hline
        $\varepsilon$ & 0.025 & 0.03316 & 0.04002 & 0.04592 & 0.05116 & 0.05592 \\
        \hline
    \end{tabular}
    \caption{Values for $\varepsilon$ for each $g_m$.}
    \label{table:eps}
\end{table}
The resulting functions of different smoothness are plotted in Figure~\ref{fig:ramps_all}.
As we can see from the left graph in Figure~\ref{fig:ramps}, the shape of all functions is quite similar.
The main difference lies in the regularity as can be seen in the zoom in Figure~\ref{fig:ramps2}.
%

\begin{figure}[!htb]
    \captionsetup[subfigure]{width=0.9\linewidth}
    \centering
    \begin{subfigure}{.5\textwidth}
        \centering
        \includegraphics[width=1\linewidth]{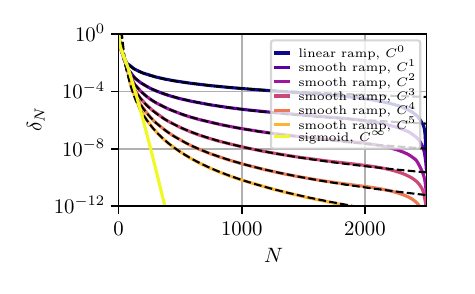}
        \caption{$\delta_N$ and error estimation for $C^{m-1}$-functions.}
        \label{fig:d_N_C_k}
    \end{subfigure}%
    \begin{subfigure}{.5\textwidth}
        \centering
        \includegraphics[width=1\linewidth]{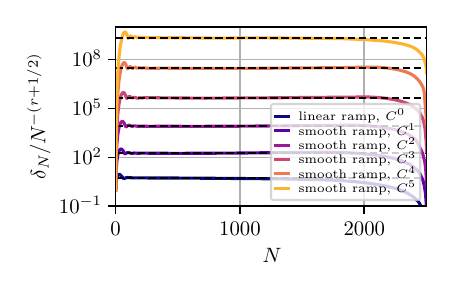}
        \caption{Estimation of $c_{r}$ s.t. $\delta_N \approx c_m N^{-(m+1/2)}$}.
        \label{fig:C_C_k}
    \end{subfigure}
    \caption{$N$-width for ramps with varying regularity.}
    \label{fig:error_bounds:2}
\end{figure}

The results concerning the $N$-width are shown in Figure~\ref{fig:error_bounds:2}.
On the left, in Figure~\ref{fig:d_N_C_k}, we compare the $N$-width $\delta_N( \cU^g )$ for $g_{q_{ m}}\in C^{m}(\Omega)$, $m=0,...,5$ and also for the $C^\infty$-sigmoid function (yellow) with exponential decay.
We also indicate the error bound from Theorem~\ref{Thm-decay-lower-bound}, i.e., $\tilde{c}_m N^{-(m+1/2)}$.
As there is no difference visible, in Figure~\ref{fig:C_C_k}, we plot the ratio of the numerically computed error and $c_m N^{-(m+1/2)}$ for a fitted $c_m$ for $m = 0,...,5$.
We see very good matches indicating that our bounds are sharp regarding $N$, in particular since the displayed functions are expected to have the Sobolev regularity $r(m)=m+1/2-\varepsilon$, see Remark \ref{rem-Hrminuseps}.
As with the decay of the jump discontinuity (c.f. Figure~\ref{fig:error_heaviside}), the numerically computed decay for $N$ close to $\min(n_x, n_\mu)$ suffers from inaccuracies that are related to the discretization error.
%

\subsection{The impact of the slope}
\begin{figure}[!htb]
    \centering
    \begin{subfigure}{.5\textwidth}\centering
        \includegraphics[width=1\linewidth]{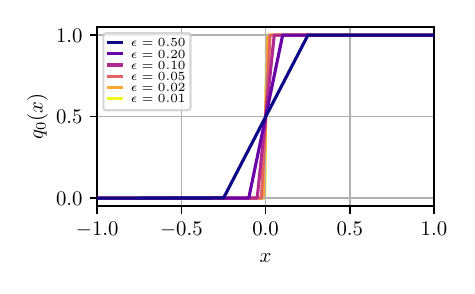}
        \caption{Function $q_0$ for different $\varepsilon$.\label{fig:steepness_f}}
    \end{subfigure}%
    \begin{subfigure}{.5\textwidth}\centering
        \includegraphics[width=1\linewidth]{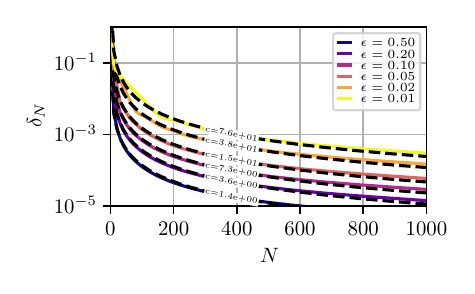}
        \caption{$\delta_N(\cU^g)$ for different $\varepsilon$.\label{fig:steepness_e}}
    \end{subfigure}
    \caption{$N$-width depending on the slope of a continuous, piecewise linear function.\label{fig:steepness}}
\end{figure}

In \S\ref{SubSec:steep} we have investigated functions with an almost identical ramp but with different smoothness.
Now, we fix the regularity and vary the slope, i.e., the maximal value of the derivative (or norm of the gradient in higher dimensions).
From our theoretical findings, we expect that the asymptotic decay rate should not be influenced by the slope.
However, all estimates involve a multiplicative factor, which might depend on the slope.
In order to clarify this, we consider a continuous, piecewise linear function with varying steepness.
We choose the function $q_0$ in \eqref{equ:polynomials:0} for different values of $\varepsilon$, see Figure~\ref{fig:steepness_f}.
The results are displayed in Figure~\ref{fig:steepness}.
We observe that the asymptotic rate is in fact identical, but the multiplicative factor grows when $\varepsilon$ decreases: the steeper the slope, the larger the $N$-widths.
%

\subsection{Beyond symmetry}\label{sec:beyhond_symm}
\begin{figure}[!htb]
    \captionsetup[subfigure]{width=0.8\linewidth}
    \begin{subfigure}[t]{.5\textwidth}\centering
        \includegraphics[width=\linewidth]{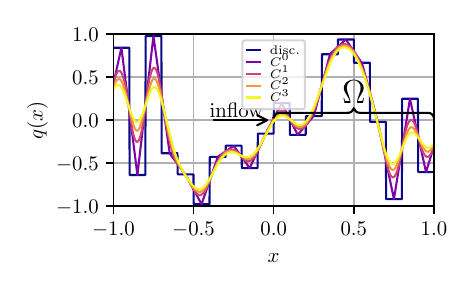}
        \caption{Piecewise constant with $20$ random steps (dark blue); increasing smoothness $C^m(\OP)$ by $m+1$-fold convolution, $m=0,...,3$.}
        \label{fig:random2_1D}
    \end{subfigure}%
    \begin{subfigure}[t]{.5\textwidth}\centering
        \includegraphics[width=\linewidth]{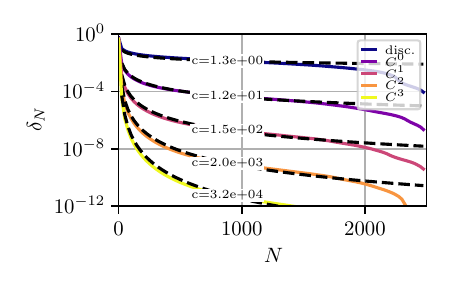}
        \caption{$N$-width decay for discontinuous function and $g\in C^m(\OP)$, $m=0,...,3$.}
        \label{fig:random_1D_decay}
    \end{subfigure}%
    \caption{$\delta_N(\cU^g)$ for random functions of different smoothness.}
    \label{fig:RandomIC}
\end{figure}

Finally, we consider almost arbitrary functions $g$ to define initial and boundary conditions for our original linear transport problem \eqref{eq:transport} in the sense that $g_{|[-1,0]}$ defines the inflow (i.e., the boundary condition) and $g_{|[0,1]}$ is the initial condition on $\Omega$.
Again, we focus on the influence of the regularity on the decay of the $N$-width.
To this end, we start by a piecewise constant discontinuous function as displayed in Figure~\ref{fig:random2_1D} (dark blue), where the height of the $20$ steps are chosen at random.
Smother versions are constructed by applying a convolution with a uniform box kernel, that is as wide as the distance between two discontinuities, see also Figure~\ref{fig:random2_1D}.\footnote{ A closer agreement between the original and the convoluted function as well as a faster error decay could be achieved through a convolution by a narrow Gaussian kernel. However, we aimed at highlighting the effect of regularity.}
The $N$-width is shown in Figure~\ref{fig:random_1D_decay}, where $-$again$-$ we clearly see the dependence of the decay on the regularity; the smoother the function, the better the rate. The rates are the same as in the previous case, but the constants (indicated in Figure~\ref{fig:random_1D_decay}) differ. This experiment confirms our results also beyond half-wave symmetry, which we had to assume for the given proofs.
%

\subsubsection*{A 2D-example}
\begin{figure}[!htb]
    \captionsetup[subfigure]{width=0.8\linewidth}
    \centering
    \begin{subfigure}{.5\textwidth}\centering
        \includegraphics[width=1\linewidth,height=38mm]{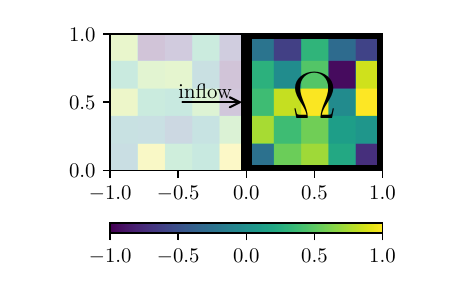}
        \caption{Piecewise constant initial and boundary conditions indicated by color boxes.\label{fig:random0_2D}}
    \end{subfigure}%
    \hfill
    \begin{subfigure}{.5\textwidth}\centering
        \includegraphics[width=1\linewidth]{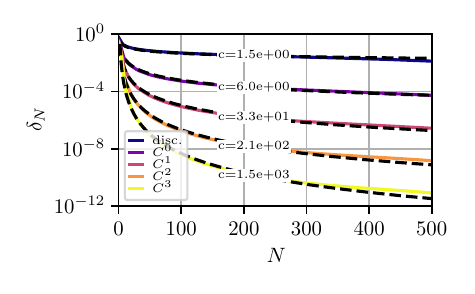}
        \caption{Kolmorogov $N$-width for $C^m(\OP)$-conditions, $m=0,...,3$.\label{fig:random_2D_decay}}
    \end{subfigure}
    \caption{2D-transport problem: $\delta_N$-width for initial- and boundary conditions of different regularity $C^m(\OP)$, $m=0,...,3$.\label{fig:random}}
\end{figure}

All our analysis above was restricted to the 1D case $\Omega=(0,1)$.
However, from the presentation it should be clear that at least some of what has been presented can be generalized to the higher-dimensional case by means of tensor products.
In order to show this also numerically, we consider a linear transport problem on $\Omega=(0,1)^2$, see Figure~\ref{fig:random0_2D}. Note, that the parameter $\mu$ remains univariate.
There, we indicate piecewise constant boundary conditions (on the left square yielding the inflow conditions) and initial conditions on $\Omega$ (on the right square).
As before, we realize initial and boundary conditions of higher regularity by applying convolutions.
The resulting $N$-widths are displayed in Figure~\ref{fig:random_2D_decay}, where we see once more that the rate is correlated to the regularity.

\section{Conclusions}\label{Sec:Outlook}
We have derived both exact representations as well as sharp bounds \rev{of} the $N$-widths for significant classes of functions \rev{$g$} used as initial and boundary values for the linear transport problem. It became clear that \rev{a parametric problem inherits its $N$-width decay by its eigenvalue decay, which equal Fourier coefficients of the data $g$ for the linear transport setting.} The influence of the regularity of $g$ on the decay has been rigorously investigated. It became clear that a poor decay of the $N$-width is \emph{only} a question of the smoothness of the solution in terms of the parameter, not of the problem itself.
In other words, the \rev{$N$-width} decay does not necessarily depend on the PDE \rev{alone}, but on the data such as initial and boundary values.
We have also seen that the constant in the decay estimate depends on the slope of the function in a severe manner.
\rev{The numerical experiments have also demonstrated, that  small changes of the data that increase regularity, can lead to vastly faster error decays.}

Our main tool is Fourier analysis and the notion of half-wave symmetric functions.
This notion allowed us to construct linear spaces \rev{which are shift-isometric spectral spaces and therefore} optimal in the sense of Kolmogorov.
Since any function can be written as a sum of even and odd HWS function, we derived a general \emph{upper} estimate for the $N$-width. 
We have investigated both the $L_2(\cP)$-based $N$-width $\delta_N(\cU^g)$ and the $L_\infty(\cP)$-based (worst case) $N$-width $d_N(\cU^g)$ and we have proven  $\delta_N(\cU^g)=d_N(\cU^g)$ for \rev{shift-isometric spaces}. 
Finally, \rev{some ideas of} the presented approach could also be generalized and adapted for other kinds of PPDEs.

\section*{Statements and Declarations}
FA gratefully acknowledges STIPINST funding [318024] from the Research Council of Norway.
The authors declare no competing interests.
\bibliography{AGU_references_v3}

\end{document}